\documentclass[]{amsart}
\usepackage[dvips]{epsfig}
\usepackage{fullpage}
\usepackage{latexsym}
\usepackage{amsmath}

\usepackage{graphicx}

\usepackage{amsmath}
\usepackage{amssymb}
\usepackage{amsbsy}

\newcommand{\vep}{\varepsilon}

\newcommand{\vex}{{\bx}}

\newcommand{\trp}{{\mbox{\scriptsize\sf T}}}
\newcommand{\vnu}{{\boldsymbol \nu}}

\newcommand{\bA}{{\boldsymbol A}}
\newcommand{\bB}{{\boldsymbol B}}
\newcommand{\bC}{{\boldsymbol C}}
\newcommand{\bD}{{\boldsymbol D}}

\newcommand{\bI}{{\boldsymbol I}}

\newcommand{\bT}{{\boldsymbol T}}

\newcommand{\bq}{{\boldsymbol q}}

\newcommand{\fbf}{{\boldsymbol f}}
\newcommand{\bg}{{\boldsymbol g}}
\newcommand{\bh}{{\boldsymbol h}}
\newcommand{\bn}{{\boldsymbol n}}
\newcommand{\bm}{{\boldsymbol m}}

\newcommand{\bu}{{\boldsymbol u}}
\newcommand{\bv}{{\boldsymbol v}}
\newcommand{\bw}{{\boldsymbol w}}
\newcommand{\bx}{{\boldsymbol x}}

\newcommand{\bz}{{\boldsymbol z}}

\newcommand{\bnu}{{\boldsymbol \nu}}

\newcommand{\bpsi}{{\boldsymbol \psi}}

\newtheorem{theorem}{Theorem}[section]
\newtheorem{lemma}[theorem]{Lemma}
\newtheorem{proposition}[theorem]{Proposition}
\theoremstyle{definition}
\newtheorem{definition}[theorem]{Definition}

\theoremstyle{remark}

\numberwithin{equation}{section}

\title{G-convergence and homogenization of viscoelastic flows.}
\author{Alexander Panchenko}
\address{Department of Mathematics, Washington State University, Pullman, WA, 99164}
\email{panchenko@math.wsu.edu}
\begin{document}
\begin{abstract}
The paper is devoted to homogenization of two-phase 
incompressible viscoelastic flows with disordered microstructure. 
We study two cases. In the first case, both phases are modeled as Kelvin-Voight viscoelastic materials.
In the second case, one phase is a Kelvin-Voight material, and the other is a viscous Newtonian fluid.
The microscale system contains the conservation of mass and balance
of momentum equations. The inertial terms in the momentum equation incorporate the actual interface
advected by the flow. In
the constitutive equations, a frozen interface is employed. The interface geometry is arbitrary: we do not assume periodicity, statistical homogeneity or scale separation. The
problem is homogenized using $G$-convergence and oscillating test functions. Since the microscale system is not parabolic, previously known constructions of the test functions do not work here. The test functions developed in the paper are non-local in time and satisfy divergence-free constraint exactly. The latter feature enables us to avoid working with pressure directly. We show that the effective medium is a single phase viscoelastic material that is not necessarily of Kelvin-Voight type. The effective constitutive equation contains a long memory viscoelastic term, as well as instantaneous elastic and viscous terms. 
\end{abstract}

\maketitle


\section{Introduction}
\noindent



Formulation of constitutive equations of multiphase materials under flow is a fundamental problem of continuum mechanics. Using mathematical homogenization theory \cite{BLP}, \cite{JKO-94}, \cite{SP-80} to solve this problem is intuitively 
appealing, but not easy. The reasons for this can
be summarized as follows: it is difficult to homogenize
evolution equations, non-linear equations, and equations involving general \footnote{not necessarily random homogeneous} geometric distribution of the constituents. In particular, addressing the last difficulty is the necessary first step in developing a complete homogenization theory for moving interface problems.

Consider a composite material with two constituents which we call phases. During flow, the interface between the phases is advected by the flow velocity. Therefore, the interface motion is coupled to the flow dynamics. A priori, one cannot expect
that a geometry that is random homogeneous at the initial time remains random homogeneous at future times.
Scale separation also cannot be expected to hold for all times in the interval of interest.
Therefore, homogenization techniques that require a specific type of geometry, e.g. two-scale convergence
(\cite{Allaire, Ngue}) and ergodic theorems, 
cannot be used in general. This leaves $G$-convergence \cite{MT, Pankov, Sp1, Sp2, jko0, jko1} and, whenever variational formulation is available,
$\Gamma$-convergence (\cite{braides, dalmaso}). The problem studied in this paper is not variational, so $G$-convergence is chosen as the main technical tool.

$G$-convergence is a general notion of functional analysis
and operator theory. In fact, $G$-convergence of an operator sequence $A^\vep$ can be identified with the 
Painleve-Kuratowski set convergence \cite{K, RW} of the corresponding operator graphs $\Gamma^\vep$
(see e.g. the definition of $G$-convergence in \cite{Pankov}). Set convergence is sequentially compact provided
the topology of $\Gamma^\vep$ has a countable base. Therefore, existence of an abstract $G$-limit operator $A$ is easily obtained for  
$A^\vep$ that may be non-linear, non-local and multi-valued. Once existence of $A$ is established, one needs to 
describe the structure of $A$, which is a problem of {\it characterization}. To solve this problem we use the method of oscillating test functions \cite{Evans1, Evans2, MT, Pankov}.
Our point of view is somewhat different from the standard one. To explain this further, consider
a sequence of problems
$$
A^\vep \bu^\vep=\fbf,
$$
where
$A^\vep: X\to X^\star$ are linear operators, $X$ is a separable Banach space and $X^\star$ is its topological dual. 
Suppose that $A^\vep$ $G$-converges to an operator $A$. This means that $\bu^\vep\to \overline{\bu}$ weakly in $X$ and $A\overline{\bu}=\fbf$
for each $\fbf\in X^\star$.
Instead of choosing a sequence $\bw^\vep$ of oscillating test functions, we prescribe a sequence of {\it corrector operators} $N^\vep: X\to X$ and then define the test functions by
$$
\bw^\vep=N^\vep\bw,
$$
where $\bw$ is an arbitrary fixed test function from a dense subset of $X$. 

It turns out that the characterization problem can be solved if $N^\vep$ satisfy\newline
\begin{equation}
\label{c1}
N^\vep\bw \to \bw~~~\mbox{weakly~in}~X,
\end{equation}
\begin{equation}
\label{c2}
\left(A^\vep\right)^\trp N^\vep \bw\to \overline{A}^\trp \bw~~~~\mbox{strongly~in}~X^\star, 
\end{equation}
possibly along a subsequence. Here, $\left(A^\vep\right)^\trp$ is the formal adjoint of $A^\vep$.
The operator $\overline{A}^\trp$ must be described as explicitly as possible. 

Indeed, if (\ref{c1}) holds then
\begin{equation}
\label{1-0}
\langle A^\vep \bu^\vep, N^\vep \bw\rangle=\langle \fbf, N^\vep \bw\rangle\to
\langle \fbf, \bw\rangle=\langle A\overline{\bu}, \bw\rangle,
\end{equation}
and if (\ref{c2}) is true,
then
\begin{equation}
\label{2-0}
\langle A^\vep \bu^\vep, N^\vep \bw\rangle=\langle \bu^\vep, \left(A^\vep\right)^\trp N^\vep \bw\rangle\to
\langle \overline{\bu},  \overline{A}^\trp \bw\rangle=\langle \overline{A}\overline{\bu}, \bw\rangle,
\end{equation}

Since $\bw$ is any function in a dense subset of $X$, (\ref{1-0}), (\ref{2-0}) imply 
$A\overline{\bu}=\overline{A}\overline{\bu}$. 

In practice, finding a corrector operator is a difficult task, and at present there are no general recipes for doing it.
Known methods of $G$-limit characterization are specific to a certain class of problems: elliptic equations in divergence form
\cite{MT}, elliptic and parabolic problems with coercivity \cite{jko0, jko1}, linear elasticity \cite{OShY}, operator
equations with monotone \cite{chiado}, and pseudo-monotone \cite{Pankov} operators. 

In most of these cases, the $G$-limit problem has the same general structure as the $\vep$-problems. For example, mixing two linear elastic materials produces an effective material that is also linear elastic. However, in
mechanics of multi-phase flows, the phases may have different physical properties and the general nature of the effective equations is often unknown. Even in the relatively simple case of dense suspensions of rigid particles in a Stokes fluid,
effective equations are subject to debate (\cite{BP, goddard, sierou-brady}). Therefore, the first goal of homogenization would be to determine the general structure of the effective equations: one-phase or two-phase, simple or not, type of memory dependence, presence of additional state variables etc. 
$G$-convergence is particularly well suited for answering questions of this kind. 

In this paper, we use $G$-convergence to study homogenization of two-phase incompressible viscoelastic flows.
The phases are modeled as Kelvin-Voight materials:
the elastic stresses satisfy Hook's law written in the spatial (Eulerian) formulation, and the viscous stresses obey Newton's law. Despite the fact that constitutive equations for the phases look similar, this model can describe the mixture of two materials, one of which is fluid-like (small elasticity), and the other is solid-like (small viscosity). 
We also study the fluid-structure interaction problem, where one phase is a Kelvin-Voight material, and the other is a viscous Newtonian fluid.
At the initial moment of time the two phases are finely mixed, and both phases occupy connected domains. The latter assumption excludes particle flows, since particles may collide, and that requires more complicated interface conditions. The density of each phase is constant, but the initial density of the mixture is highly oscillatory. This means that the sequence of initial densities does not converge strongly. No further assumptions are made about the initial interface geometry. In particular, we do not assume periodicity, statistical homogeneity or scale separation.

The system of equations contains the mass balance and the momentum balance equations. The mass balance is needed because the initial density of the mixture is not constant, and the interface geometry changes in time.
We also make the following choice regarding modeling of the interface. The mass balance equation and
the convective terms in the momentum equation incorporate a moving interface advected by the flow velocity. 
The constitutive equations for the stress employ a frozen interface, that is the interface that existed at the initial time. 
This choice makes sense both physically and mathematically.

From the point of view of physics, freezing the interface in constitutive equations is the only option compatible with 
the Hook's law. In linear elasticity, deformations are assumed to be small, and thus spatial and referential descriptions are identified. In the referential description, the interface is always fixed. The reason for this is simple: a material particle that belongs initially to phase one remains in that phase at future times. Therefore, in the framework of linear elasticity, the interface must be fixed in both descriptions. Otherwise elastic forces may lead to non-physical energy dissipation. More details illustrating this point are given in Appendix A. 

This may seem like a reason to use fixed interface in all terms and work with the
formally linearized equations of motion. This is not satisfactory for viscoelastic composites,
because convective effects may be significant even when the interface is nearly stationary. Indeed, when
one phase is fluid-like, it may flow at large Reynolds number in the cavities of a matrix made from the second, solid-like phase.
In that case the deformations inside the solid-like phase are small, as are the deformations of the interface. However, the deformations in the fluid-like phase may be large, and convective effects inside that phase cannot be neglected.
Formal calculations, presented in Appendix
\ref{ap-model}, show that presence of convective terms is incompatible with the frozen
interface. It is also impossible to meet the natural interface conditions unless
convective terms are present in both phases, and not just in a fluid-like phase.

To summarize, the microscale model in the paper partially accounts
for the shape variation of the interface and satisfies all of the above natural requirements. In addition, these equations provide a good benchmark model
for studying more complicated problems with moving interface such as flows
of two immiscible fluids with surface tension, flows of nonlinear viscoelastic materials, fluid-particle flows etc.

Mathematical study of the above microscale problem makes sense because little is known about $G$-limits of non-parabolic systems of fluid mechanics. Since the sequences of initial data do not converge strongly, the known compactness results for 
weak solutions such as Theorem 2. 4.
in \cite{lions2} and Theorem 5.1 in \cite{lions3} cannot be expected to hold. In a sense, the problems of the type studied here are more difficult than existence questions. Because of the lack of compactness, 
the structure of $G$-limits for generic continuum mechanics systems may be different from the structure of $\vep$-problems, but this is not known. We show that in the present case, a new long memory term appears in
the effective constitutive equations. Such results were previously obtained for linear equations with fixed interface. Formally this was done for locally periodic geometry in \cite{SP-80} and \cite{BK-81}, and rigorously for periodic geometry in \cite{GM-00} using two-scale convergence. In \cite{GP-04}, oscillating test functions were used for non-periodic scale separated geometries (some details of the proofs in that paper were omitted). In \cite{GM-00} and \cite{GP-04}, the oscillating test functions did not satisfy divergence-free constraint. 
Using arbitrary test functions in incompressible problems makes it necessary to work with pressure directly. This is a serious obstacle, since good estimates on the pressure are not available even for density dependent Navier Stokes equations \cite{lions2}. 
In this paper, we incorporate physically meaningful convective terms, allow for interfacial motions, and propose a construction of divergence-free oscillating test functions. The construction of the corrector operators involves certain auxiliary functions satisfying auxiliary problems. The right hand sides in auxiliary problems are chosen to satisfy condition (\ref{c1}) on the corrector operators $N^\vep$. The treatment here is inspired by the "condition N" in the papers \cite{jko0, jko1} of Zhikov, Kozlov and Oleinik (see also the book \cite{OShY} for an application to linear elasticity). It seems that many ideas in these works can be extended to non-parabolic evolution problems, including problems with moving interface.

The main result of the paper are Theorem \ref{main-thm} in Sect. \ref{sect:outline}, and Theorem \ref{fsi-main} in 
Sect. \ref{sect:fsi}.
There we show that the effective system of equations describes a single-phase
incompressible viscoelastic material. The effective system contains
equations of mass and momentum balance. The effective constitutive equations contain a linearly elastic term, a linearly viscous term and a viscoelastic term that models a long term memory dependence of the effective material.

The paper is organized as follows. Section 1 is introductory. Microscale problem formulation and properties of finite energy weak solutions are given in Section 2. An outline of the existence proof for each fixed $\vep>0$ is presented in Appendix \ref{ap-exist}. With some minor changes due to the presence of the elastic stress,
we follow closely the existence proof for incompressible density-dependent
Navier-Stokes equations in \cite{lions2}, Sect. 2.3, 2.4. The global 
weak solutions of Leray type
obtained in this way satisfy the energy inequality. Section 3-6 contain a detailed study of the the case of two Kevin-Voight materials. In Section 3, we formulate the main theorem and provide an outline of the proof. Section 4 is devoted to constructing the corrector operators and oscillating test functions. Section 5 deals with passage to the limit in the inertial terms of the momentum balance equation. In Section 6 we obtain the effective constitutive equations and combine all the results to finish the proof of the main theorem. Finally, in Section 7 we indicate the changes necessary to treat the fluid-structure interaction case, and state the main theorem for this case.

\section{Micro-scale problem.}

\subsection{Equations of balance and constitutive equations}
\subsubsection{Choice of a model}
We consider two-phase materials in which at least one of the phases resist
shearing, and the material stress tensor can be written as a sum of an elastic
(conservative) and dissipative stresses. To avoid mathematical difficulties
in dealing with nonlinear elasticity we limit our investigation to flows
for which the Hook's law of linear elasticity is an appropriate model of the
elastic stress. We further suppose that deformation of the interface is small.
In this case, the equations of motion are often formally linearized, assuming that the density is constant, and spatial and referential description are identified. As a consequence, the interface (always fixed in the referential description) is fixed in the spatial description. This is unsatisfactory when the physical properties of the phases are different., e.g. one phase is solid-like and the other one fluid-like. In this case, large deformations of the fluid-like phase may occur even when the deformations in the solid-like phase are small. Consequently, the contribution of the inertial terms to the overall momentum balance cannot be neglected. In addition, a correct model should respect physically realistic interface conditions (continuity of velocity and equality of tractions) as well as physically correct
(formal) energy balance. These requirements make it necessary to
consider inertia in both phases, with the corresponding terms that include
a moving interface and densities governed by the mass conservation equations. 
By contrast, Hook's law for the elastic stress may lead to the 
non-physical dissipation of the elastic energy unless the interface in the constitutive equations is frozen. In Appendix
A we present some (formal) calculations illustrating the above points.

In this paper, we work with the micro-scale equations that employ a moving interface
in the inertial terms and a frozen interface in the constitutive equations
for the stress. This model can be viewed as a transition between a completely
linear model with constant densities and static interface, and a fully
nonlinear nonlinear that involves nonlinear constitutive equations for
the elastic stress.

\subsubsection{Micro-scale equations}
Let $\rho^\vep, {\bv}^\vep$ denote, respectively, the density and velocity of the composite.
We also define 
\begin{equation}
\label{displ}
{\bu}^\vep(t, {\bx})=\int_0^t {\bv}^\vep(\tau, {\bx})d\tau.
\end{equation}

\noindent
{\bf Interface evolution equation.}
Evolution of $V^\vep(t)$, $W^\vep(t)$
can be described by the interface evolution equation. Let $\theta^\vep(t,{\bx})$ denote the
characteristic function of $V^\vep(t)$. In the referential description, $V^\vep$ is fixed. Therefore, the
material derivative of $\theta^\vep$ is zero. In the spatial description we have
\begin{equation}
\label{int-cons}
\partial_t \theta^\vep+{\bv}^\vep\cdot \nabla \theta^\vep=0~~~{\rm in}~U.
\end{equation}
Equation (\ref{int-cons}) is supplemented with the initial condition
\begin{equation}
\label{in-theta}
\theta^\vep(0, {\bx})=\theta_0^\vep(\bx)~~~{\rm in}~U.
\end{equation}

\noindent
{\bf Mass conservation equation.}
Next, we state the mass conservation equation. 
The composite density satisfies 
\begin{equation}
\label{t-m}
\partial_t \rho^\vep+{\rm div}(\rho^\vep{\bv}^\vep)=0,~~~{\rm in}~U
\end{equation}
with the initial condition
\begin{equation}
\label{t-m-in}
\rho^\vep(0, {\bx})=\rho_1\theta_0^\vep(\bx)+
\rho_2(1-\theta_0^\vep(\bx)),~~~{\bx}\in U,
\end{equation}
where $\rho_1$, $\rho_2$ are the densities of
the respective phases. These densities are assumed to be constant and bounded below by a positive constant.

\noindent
{\bf Incompressibility:} 
\begin{equation}
\label{0div}
{\rm div}~{\bv}^\vep=0.
\end{equation}

\noindent
{\bf Momentum balance for the composite:} 
\begin{equation}
\label{mb}
\partial_t {\bv}^\vep+{\rm div}(\rho^\vep {\bv}^\vep \otimes {\bv}^\vep)-{\rm div}\left(
{\boldsymbol T}_1^\vep-P_1^\vep\boldsymbol I+{\boldsymbol T}_2^\vep-P_2^\vep{\boldsymbol I}\right)=0.
\end{equation}
Here ${\boldsymbol T}_s^\vep-P_s^\vep{\boldsymbol I}$ is the stress tensor in the phase $s$, 
$P_s^\vep$ is the pressure, and ${\boldsymbol I}$ denotes the unit tensor. The
initial conditions for (\ref{mb}) are
\begin{equation}
\label{mb-ic}
{\bv}^\vep(0, {\bx})={\bv}_0({\bx}),
\end{equation}
where ${\bv}_0$ does not depend on $\vep$.
In addition, (\ref{displ}) implies that ${\bu}^\vep(0, {\bx})=0$.

On the boundary $\partial U$, the condition 
\begin{equation}
\label{zero-v}
{\bv}^\vep(t, {\bx})=0
\end{equation} 
is imposed.

\noindent
{\bf Constitutive equations.} 
As explained above, static interface seems to be a natural choice that is relatively easy to handle (unlike combining
referential formulation for the elastic stress with the spatial formulation
for the viscous stress), and 
compatible with the spatial form of the Hook's law for the elastic part of the stress. 
We therefore define
\begin{eqnarray}
&&\boldsymbol T^\vep={\boldsymbol T}_1^\vep+{\boldsymbol T}_2^\vep-P^\vep{\boldsymbol I}, ~{\rm where}~ P^\vep=P_1^\vep+P_2^\vep, \label{ac-const}\\
&&{\boldsymbol T}_1^\vep(t,{\bx})=\theta_0^\vep(\bx)\left({\boldsymbol A}^1e({\bu}^\vep)+{\boldsymbol B}^1e({\bv}^\vep)\right), \nonumber\\
&&{\boldsymbol T}_2^\vep(t,{\bx})=(1-\theta_0^\vep(\bx))\left({\boldsymbol A}^2e({\bu}^\vep)+{\boldsymbol B}^2e({\bv}^\vep)\right),\nonumber
\end{eqnarray}
and ${\boldsymbol A}^s, {\boldsymbol B}^s$, $s=1,2$ are constant material tensors. In (\ref{ac-const}), $e=\frac 12(\nabla+\nabla^T)$ is the
symmetric part of the gradient. We assume that 
both phases are isotropic. In that case
\begin{equation}
\label{material}
A^s_{ijkl}=\mu^s \delta_{ik}\delta_{jl},
~~~~~
B^s_{ijkl}=\nu^s \delta_{ik}\delta_{jl},
~~~~s=1,2,
\end{equation}
where $\mu^s$ are the elastic moduli and $\nu^s$ are the viscosities of the phases. All these constants are supposed to be positive.

We also assume that the tensors $\boldsymbol A^\vep, \boldsymbol B^\vep$ satisfy
\begin{equation}
\label{ellipticity}
\alpha_1 \boldsymbol \xi\cdot \boldsymbol \xi\leq \boldsymbol A^\vep\boldsymbol \xi\cdot \boldsymbol \xi\leq \alpha_2 \boldsymbol \xi\cdot \boldsymbol \xi,
~~~~~~~~~~
\beta_1 \boldsymbol \xi\cdot \boldsymbol \xi\leq \boldsymbol B^\vep\boldsymbol \xi\cdot \boldsymbol \xi\leq \beta_2 \boldsymbol \xi\cdot \boldsymbol \xi.
\end{equation}
for each $\boldsymbol \xi\in {\bf R}^{3\times 3}$, with $\alpha_i>0, \beta_i>0$ independent of $\vep$.

\noindent
{\it Remark}. Together, (\ref{ac-const}), 
(\ref{t-m}), (\ref{mb}) form a closed system, and thus the 
interface evolution equation (\ref{int-cons}) can be dropped. This fact is important for compressible flows for which the mass conservation equation is stable with respect to weak convergence, while equation (\ref{int-cons}) is not.
For incompressible flows considered in this paper, the interface evolution equation has the same structure $\partial_t+{\rm div}({\bv}^\vep \cdot)$ as the mass conservation. Moreover, if the initial densities $\rho_1, \rho_2$ are constant, then the densities of the phases remain constant during the motion. In this case, the interface evolution equation and the mass conservation equation are essentially equivalent.

\noindent
{\bf Interface conditions}. There are two  
interface conditions: the first is continuity of ${\bv}$ across the interface (which
is the actual moving interface governed by (\ref{int-cons})), and the second 
is the equality of tractions $\left( {\boldsymbol T}^\vep_s-P_s^\vep\bI\right) \vnu_s$ on the frozen interface. 
Here $\vnu_s$ denotes the exterior (to the phase $s$) unit normal to the frozen interface.

\subsubsection{Weak formulation of the micro-scale problem.}
In this section we provide the weak formulation of the problem to be homogenized.
It consists of the mass conservation and momentum balance equations.\newline
\noindent 
{\bf Mass conservation.} 
\begin{equation}
\label{1.2}
\int_U \rho^\vep(0, \vex)\phi(0,\vex)d\vex-\int_{I_T}\int_U \rho^\vep \partial_t \phi d\vex dt-
\int_{I_T}\int_U \rho^\vep {\bv}^\vep\cdot\nabla \phi d\vex dt=0,
\end{equation}
where $\rho^\vep(0, \vex)$ is given by (\ref{t-m-in}), and
$$
I_T=[0, T).
$$
Equation (\ref{1.2}) is supposed to hold for each smooth test function 
$\phi$, equal to zero on $\partial U$ and vanishing for $t\geq T$.  

\noindent
{\bf Momentum balance.}
\begin{eqnarray}
&&-\int_U \rho^\vep(0, \vex){\bv}_0 \cdot {\boldsymbol \psi} d\vex-\int_{I_T}\int_U \rho^\vep 
{\bv}^\vep\cdot\partial_t {\boldsymbol
\psi} d\vex dt
-\int_{I_T}\int_U \rho^\vep {\bv}^\vep\otimes {\bv}^\vep\cdot \nabla{\boldsymbol \psi}d\vex dt
\label{1.3}
\\
&&+\int_{I_T}\int_U ({\boldsymbol T}_1^\vep+{\boldsymbol T}_2^\vep)\cdot e({\boldsymbol \psi})d\vex dt
=0.\nonumber
\end{eqnarray}
Equation (\ref{1.3}) holds
for all smooth test-functions ${\boldsymbol \psi}$, such that ${\rm div}\;\bpsi=0$, $\bpsi$ equal to zero
on $\partial U$, ${\boldsymbol \psi}(t, {\bx})=0$ when $t\geq T$. The dependence of $\bT_i^\vep$ on
$\bv_i^\vep$ and $\theta^\vep_0$ is given by (\ref{ac-const}).

\noindent
{\it Remark}. It is important to note that, because of the condition $\bv^\vep=0$ on $\partial U$, the identity (\ref{1.3}) also holds for test functions $\bpsi$ with the condition $\bpsi=0$ replaced by a less restrictive $\bpsi\cdot \bnu=0$ on $\partial U$, $\bnu$ is the exterior unit normal to $\partial U$. This fact will be used in Section  \ref{sect:corrector} to construct oscillating test functions.

\subsection{Finite energy weak solutions and bounds}
\label{weak-sect}
We suppose that the initial conditions satisfy
\begin{eqnarray}
&& 0<C_1\leq \rho^\vep(0, {\bx})\leq C_2,\label{in-den}\\
&& {\bv}_0\in H_0^1(U)\label{in-vel}
\end{eqnarray}
with $C_1, C_2$ independent of $\vep$. 
The system (\ref{1.2}), (\ref{1.3}) closely resembles the system of density-dependent Navier-Stokes equations with density-dependent viscosity studied
in \cite{lions2}, ch. 2. The only difference is the presence of the strain-dependent terms of the type ${\boldsymbol A}e({\bu}^\vep)$ in the constitutive
equations. When the viscosity does not vanish, as is the case here, the existence
and overall properties of the weak solutions are determined by the viscosity, and not elasticity, of the medium. The proof of existence for each fixed $\vep>0$ is outlined
in Appendix \ref{ap-exist}. 
It yields, 
for each $\vep\in \{\vep_k\}_{k=1}^\infty$, existence of the finite energy weak solutions of (\ref{1.2}), (\ref{1.3}) with the same properties
as in \cite{lions2}, theorem 2.1, namely
\begin{eqnarray}
&& \rho^\vep \in L^\infty(I_T\times U)\cap C(I_T, L^p(U)),~~~{\rm for}~{\rm all}~
1\leq p<\infty,
\label{den}\\
&& {\bv}^\vep \in L^2(I_T, H_0^1(U)),\label{vel}
\end{eqnarray}
satisfying the energy inequalities
\begin{eqnarray}
&&\frac 12\int_U \rho^\vep|{\bv}^\vep|^2 d{\bx}(t)+
\int_U \left[{\boldsymbol A}^1\theta_0^\vep(\bx)
+{\boldsymbol A}^2\left(1-\theta_0^\vep(\bx)\right)\right]
|e({\bu}^\vep)|^2 d{\bx}(t)\label{energy1}\\
&&+\int_{I_T}\int_U \left[{\boldsymbol B}^1\theta_0^\vep(\bx)
+{\boldsymbol B}^2\left(1-\theta_0^\vep(\bx)\right)\right]
|e({\bv}^\vep)|^2 d{\bx}dt\leq \nonumber\\
&& \frac 12\int_U \rho^\vep|{\bv}^\vep|^2 d{\bx}(0)+
\int_U \left[{\boldsymbol A}^1\theta_0^\vep(\bx)
+{\boldsymbol A}^2\left(1-\theta_0^\vep(\bx)\right)\right]
|e({\bu}^\vep)|^2 d{\bx}(0).
\nonumber
\end {eqnarray}

Let us list the implications of the above estimates. First, renormalizing solutions of the
mass conservation equations, as in \cite{lions2}, sect 2.3, we obtain
\begin{equation}
\label{ap-den}
\parallel \rho^\vep\parallel_{L^\infty(I_T\times U)}\leq C_2,
\end{equation}
Next we note that (\ref{displ}) implies ${\bu}^\vep(0, {\bx})=0$, and by
(\ref{in-den}), (\ref{in-vel}), the other initial conditions are bounded independent of $\vep$. This implies that the left hand side of (\ref{energy1}) is bounded independent of $\vep$, so that
\begin{equation}
\label{apriori-stress}
\parallel{\boldsymbol T}^\vep\parallel_{L^2(I_T\times U)}\leq C,
\end{equation}
with $C$ independent of $\vep$.
Combining (\ref{apriori-stress}) with the first Korn inequality for
functions with zero trace on the boundary (see, e.g. \cite{OShY}, th. 2.1), and then with Poincar\'e inequality, we deduce
\begin{equation}
\label{ap-vel}
\parallel {\bv}^\vep\parallel_{L^2(I_T, H_0^1(U))}\leq C,
\end{equation} 
\begin{equation}
\label{ap-dis}
\parallel {\bu}^\vep\parallel_{L^\infty(I_T, H_0^1(U))}\leq C.
\end{equation} 
with $C$ independent of $\vep$. Then it follows that
\begin{equation}
\label{ap-flux}
\parallel \rho^\vep|{\bv}^\vep|^2\parallel_{L^\infty(I_T, L^1(U))}\leq C,
\end{equation}
which, together with (\ref{ap-den}) yields
\begin{equation}
\label{ap-flux1}
\parallel \rho^\vep{\bv}^\vep\parallel_{L^2(I_T, L^2(U))}\leq C,
\end{equation}
with $C$ independent of $\vep$. 

The above uniform bounds allow us to extract  a
subsequence, still denoted by $\{\vep_k\}_{k=1}^\infty$ such that
\begin{eqnarray}
&& \rho^\vep(0, {\bx})\to \overline{\rho}_0,~~~{\rm weak}-*~{\rm in}~ 
L^\infty(U),\label{lim-in-den}\\
&& \rho^\vep\to \overline{\rho}, ~~~{\rm weak}-*~{\rm in}~ L^\infty(I_T\times U),\label{lim-den}\\
&& {\bv}^\vep \to \overline{\bv},~~~{\rm weak}~{\rm in}~ L^2(I_T, H_0^1(U)),
\label{lim-vel}\\
&& {\boldsymbol T}^\vep\to \overline{\boldsymbol T},
~~~{\rm weak}~{\rm in}~ L^2(I_T,L^2(U)),\label{lim-stress}\\
&& {\bu}^\vep \to \overline{\bu}, ~~~{\rm weak}~{\rm in}~L^2(I_T, H_0^1(U)),
\label{lim-displ}\\
&& \rho^\vep{\bv}^\vep\to \overline{\bz},
~~~{\rm weak}~{\rm in}~ L^2(I_T,L^2(U)),\label{lim-flux}\\
&& P^\vep \to \overline P, ~~~{\rm weak}~{\rm in}~ L^2(I_T,L^2(U)),\label{lim-p}
\end{eqnarray}

The last comment concerns $\partial_t \rho^\vep$ and 
$\partial_t(\rho^\vep{\bv}^\vep)$. For technical reasons we need
bounds on these sequences in the space $L^1(I_T, W^{-m,1}(U))$, where
$m>0$ may be large. Such bounds can be also deduced from (\ref{energy1}) and
(\ref{1.2}), (\ref{1.3}) (see \cite{lions3}).


\section{Main theorem and outline of the proof}
\label{sect:outline}
The main theorem is as follows. 
\begin{theorem}
\label{main-thm}
The limits
$\overline{\rho}, \overline{\bv}$, $\overline{\bu}$ satisfy 
\begin{equation}
\label{main1}
{\rm div}\;\overline{\bv}=0,
\end{equation}
and the integral identities
\begin{equation}
\label{main2}
\int_U \overline{\rho}_0\phi(0,\vex)d\vex-\int_{I_T}\int_U \overline{\rho}\partial_t \phi d\vex dt-
\int_{I_T}\int_U \overline{\rho}\;\overline{\bv}\cdot\nabla \phi d\vex dt=0,
\end{equation}
\begin{eqnarray}
&&-\int_U \overline{\rho}_0{\bv}_0 \cdot {\boldsymbol \psi} d\vex-\int_{I_T}\int_U \overline{\rho}\;
\overline{\bv}\cdot\partial_t {\boldsymbol
\psi} d\vex dt
-\int_{I_T}\int_U \overline{\rho}\;\overline{\bv}\otimes \overline{\bv}\cdot \nabla{\boldsymbol \psi}d\vex dt
\label{main3}
\\
&&+\int_{I_T}\int_U \overline{\bT}\cdot e({\boldsymbol \psi})d\vex dt
=0.\nonumber
\end{eqnarray}
for all smooth test functions 
$\phi$, ${\boldsymbol \psi}$, such that ${\rm div}\;\bpsi=0$, and $\phi, \bpsi$ are equal to zero
on $\partial U$ and vanish for $t\geq T$. 

Moreover, there exist the effective tensors
$\overline{\bA}\in L^2(U), \overline{\bB}\in L^2(U)$ and $\overline{\bC}\in L^2(I_T\times U)$ such that the effective deviatoric stress $\overline{\bT}$ satisfies
 
\begin{equation}
\label{main4}
\overline{\bT}=
\overline{\bA}e(\overline{\bu})+
\overline{\bB}e(\overline{\bv})+
\int_0^t \overline{\bC}(t-\tau)e(\overline{\bv})(\tau)d\tau.
\end{equation}
\end{theorem}
\noindent
{\it Remark}. Equations (\ref{main1})--(\ref{main4}) essentially mean that the effective equations for
$\overline{\rho}, \overline{\bv}$, $\overline{\bu}, \overline{P}$
are
\begin{eqnarray}
&\partial_t \overline{\rho}+{\rm div}(\overline{\rho}\;\overline{\bv})=0,&\nonumber\\
&{\rm div}~\overline{\bv}=0,\nonumber& \\
& \partial_t(\overline{\rho}\;\overline{\bv})+
{\rm div}(\overline{\rho}\;\overline{\bv}\otimes \overline{\bv})-
{\rm div}~\overline{\boldsymbol T}+\nabla \overline{P}=0,&
\label{ef-mom}
\end{eqnarray}
with initial and boundary conditions
$$
\overline{\bv}(0,{\bx})={\bv}_0,~~~\overline{\rho}(0, {\bx})=\overline{\rho}_0,~~~
\overline{\bv}=0~~\mbox{on}~\partial U,
$$
and $\overline{\bT}$ given by (\ref{main4}).

The result follows from a number of propositions and theorems. Here, the outline of the proof is presented for the reader's convenience.

\noindent
1. {\it Construct the corrector operators $N^\vep\bw={\bw}^\vep$}.\newline
\noindent 
$N^\vep$ are defined by (\ref{appr-otf})--(\ref{phi1}). The proposed construction is non-local in $t$ and satisfies 
the divergence-free constraint. 

\noindent
2. {\it Auxiliary problems}. 
Equation (\ref{appr-otf}) includes two types of auxiliary functions: $\bn^{pq, \vep}$ and $\bm^{pq, \vep}$.
An expression of the time derivative of $\bw^\vep$ additionally contains the final value $\bm_T^{pq, \vep}=\bm^{pq, \vep}(T)$. These three families of functions satisfy auxiliary problems.
The right hand sides
in the auxiliary problems are chosen so that the sequences of solutions converge to zero weakly in the appropriate 
Sobolev type spaces. The choice of the right hand sides involves abstract $G$-limit operators corresponding to each of the auxiliary problems. An analysis of auxiliary problems is presented in Sections \ref{sect:a-prob1}, \ref{sect:a-prob2}. 
Propositions \ref{rhs-prop}, \ref{ab-rhs-choice} are used repeatedly in the remainder of the proof.

\noindent
3. {\it Convergence of $\bw^\vep$}. Using estimates from step 2, show that 
$$
\bw^\vep \to \bw,\;\;\;\; \partial_t\bw^\vep\to \partial_t\bw\;\mbox{weakly~in}~L^2(I_T, H_0^1(U)).
$$
This is done in Proposition \ref{otf-convergence}. This convergence corresponds to the condition (\ref{c1}) on $N^\vep$.

\noindent
4. {\it Convergence of the inertial terms}. Proposition \ref{otf-convergence} combined with
the Lemma \ref{lemma:lions} (Lemma 5.1 from \cite{lions2}) implies that
$$
\int_U \rho^\vep(0, {\bx}){\bv}_0\cdot{\bw}^\vep+
\int_0^T\int_U \rho^\vep{\bv}^\vep\cdot{\bw}^\vep_t+
\int_0^T\int_U \rho^\vep{\bv}^\vep\otimes {\bv}^\vep\cdot\nabla {\bw}^\vep
$$
converges to
$$
\int_U \overline\rho(0, {\bx}){\bv}_0\cdot{\bw}+
\int_0^T\int_U \overline{\rho}~\overline{\bv}\cdot{\bw}_t+
\int_0^T\int_U \overline{\rho}~\overline{\bv}\otimes \overline{\bv}
\cdot\nabla {\bw}.
$$
This is done in Proposition \ref{first-in}.

\noindent
5. {\it Compensated compactness of the stress}. Convergence of the inertial terms implies 
$$
\int_0^T\int_U \boldsymbol{T}^\vep\cdot e({\bw}^\vep)\to
\int_0^T\int_U \overline{\boldsymbol{T}}\cdot e({\bw}).
$$ 
This is shown in Proposition \ref{prop-comp}.

\noindent
6. {\it Effective stress}. Characterization of $\overline{\bT}$ is obtained in Theorem \ref{thm:effective}.
This step corresponds to establishing condition (\ref{c2}) on $N^\vep$. The main idea is to write
\begin{equation}
\label{master}
\int_0^T\int_U \boldsymbol{T}^\vep\cdot e({\bw}^\vep)=
-\langle \bu^\vep, {\rm div}\left(\bA^\vep-\bB^\vep\partial_t\right)\bw^\vep\rangle.
\end{equation}
The expression for $\langle \bu^\vep, {\rm div}\left(\bA^\vep-\bB^\vep\partial_t\right)
\bw^\vep\rangle \equiv \langle \bu^\vep, {\rm div}\left(\bA^\vep-\bB^\vep\partial_t\right)N^\vep\bw\rangle$ contains a number of terms depending on $\bn^{pq, \vep}, \bm_T^{pq, \vep}$, $\bm^{pq, \vep}$ and the corresponding pressures. In some of these terms we can pass to the limit using Lemma 5.1 from \cite{lions2}, since $\bu^\vep$ has Sobolev regularity. In other terms this is not possible, but these terms vanish by design of the auxiliary problems. 
The effective tensors 
$\overline{\bA}, \overline{\bB}$  and $\overline{\bC}$ are obtained as weak limits of the three fluxes that appear in the auxiliary problems for, respectively,
$\bn^{pq, \vep}, \bm_T^{pq, \vep}$ and $\bm^{pq, \vep}$.

\noindent
7. {\it Mass conservation}. Mass conservation equation is weakly stable. This is a well known fact (see \cite{lions2}).

\section{Corrector operators and oscillating test functions}
\label{sect:corrector}
We look for corrector operators of the form
\begin{eqnarray}
\label{appr-otf}
N^\vep\bw\equiv {\bw}^\vep(t, {\bx})& = & {\bw}(t, {\bx})+
{\bn}^{pq, \vep}({\bx})e({\bw})_{pq}\\
&&+\int_t^T \bm^{pq, \vep}(t-\tau+T, \bx)e({\bw}_t)_{pq}(\tau, {\bx})d\tau+
\nabla \phi^\vep,
\nonumber
\end{eqnarray}
(Summation over $p, q\in \{1, 2, 3\}$ is assumed).
Here, $\bw\in C_0^\infty(I_T\times U), {\rm div}\;\bw=0$ is an arbitrary test function,
$\bn^{pq, \vep}\in H_0^1(U), \bm^{pq, \vep}\in L^2(I_T, H_0^1(U))$ are to be specified. So far we require
that
\begin{equation}
\label{nodiv}
{\rm div}\;\bn^{pq, \vep}=0, ~~~~~~{\rm div}\;\bm^{pq, \vep}=0.
\end{equation}
More conditions will be imposed below. 

The function $\phi^\vep\in L^2(I_T, H^1(U))$ satisfies
\begin{eqnarray}
\label{phi1}
&& \Delta \phi^\vep=-{\bn}^{pq, \vep}\cdot \nabla e({\bw})_{pq}-
\int_t^T \bm^{pq, \vep}(t-\tau+T, \bx)\cdot \nabla e({\bw}_t)_{pq}(\tau, {\bx})d\tau,\\
&& \nabla\phi^\vep\cdot\bnu=0~{\rm on}\;\partial U.\nonumber
\end{eqnarray}

The choice of the first three terms in (\ref{appr-otf}) is motivated by similar expressions
used in periodic \cite{GM-00} and scale-separated \cite{GP-04} homogenization. The last term
are need to enforce divergence-free constraint. This is necessary in order to avoid dealing with pressure
in (\ref{1.3}) which is not $L^1_{loc}(U)$ in general.

Note also that $\nabla\phi^\vep\cdot\bnu$ is zero on $\partial U$ for almost all $t$. This makes $\bw^\vep$ correctly defined test functions for (\ref{1.3}) (see the Remark following that equation). 
Moreover, 
\begin{lemma}
\label{lemma-nodiv}
The function $\bw^\vep$ defined by (\ref{appr-otf})--(\ref{phi1}) satisfies
$$
{\rm div}\;\bw^\vep=0.
$$
\end{lemma}

\noindent
{\it Proof}. Taking divergence of (\ref{appr-otf}) and using (\ref{nodiv}) we find
\begin{eqnarray*}
&& {\rm div}\left({\bw}(t, {\bx})+{\bn}^{pq, \vep}e({\bw})_{pq}+
\int_t^T \bm^{pq, \vep}(t-\tau+T, \bx)e({\bw}_t)_{pq}(\tau, {\bx})d\tau\right)
\\
&& =
{\bn}^{pq,\vep}\cdot \nabla e({\bw})_{pq}+\int_t^T \bm^{pq, \vep}(t-\tau+T, \bx)\cdot \nabla e({\bw}_t)_{pq}(\tau, {\bx})d\tau,
\end{eqnarray*}
and the claim follows from the condition (\ref{phi1}).

\hspace{14cm}$\blacksquare$

\subsection{Auxiliary problem for $\bm^{pq, \vep}$.}
\label{sect:a-prob1}
In this subsection, $p, q$ are fixed, so we drop them to simplify notations, and write ${\bm}^\vep$ instead
of ${\bm}^{pq, \vep}$ and so on.
We look for $\bm^\vep$ that solve the auxiliary problem
\begin{equation}
\label{aux1}
-{\rm div}\left( \boldsymbol A^\vep e({\bm}^\vep)-\boldsymbol B^\vep e({\bm}_t^\vep)\right)-\nabla P^\vep_3={\fbf},~~~~~~~{\rm div}\;{\bm}^\vep=0,
\end{equation}
satisfying the condition
\begin{equation}
\label{aux2}
{\bm}^\vep(T, {\bx})={\bm}_T^\vep.
\end{equation}
The objective of this section is to show that the right hand side $\fbf$ can be chosen so that
${\bm}^\vep$ converges weakly to zero in an appropriate space.

Let $\psi(t)\in C^\infty(I_T)$ satisfy $\psi(0)=0$ and $\psi(T)=1$. We use this function to reduce (\ref{aux1}), (\ref{aux2}) to a problem with a different right hand side and zero condition at $t=T$. Writing
$$
{\bm}^\vep=\hat{\bm}^\vep+\psi{\bm}_T^\vep
$$
we deduce that $\hat{\bm}^\vep$ solves
\begin{equation}
\label{aux3}
-{\rm div}\left( \boldsymbol A^\vep e(\hat{\bm}^\vep)-\boldsymbol B^\vep e(\hat{\bm}_t^\vep)\right)-\nabla P^\vep={\fbf}+{\rm div}\left( \boldsymbol A^\vep e(\psi{\bm}^\vep_T)-
\boldsymbol B^\vep e(\psi_t{\bm}_T^\vep)\right)
,~~~~~~~{\rm div}\;\hat{\bm}^\vep=0,
\end{equation}
with the inital condition
\begin{equation}
\label{aux4}
\hat{\bm}(T, {\bx})=0.
\end{equation}

Define the spaces
\begin{equation}
\label{V}
{\mathcal V}=\{{\bv}\in L^2(I_T, H_0^1(U)), ~{\rm div}\;{\bv}=0\},
\end{equation}

\begin{equation}
\label{W}
{\mathcal W}=\{{\bv}\in {\mathcal V}: {\bv}_t\in {\mathcal V}\},
\end{equation}
\begin{equation}
\label{WT}
{\mathcal W}_T=\{{\bv}\in {\mathcal W}: {\bv}(T)=0\},
\end{equation}

The space ${\mathcal V}$ is equipped with a norm
\begin{equation}
\label{norm1}
\parallel {\bv} \parallel_{{\mathcal V}}=\left( \int_0^T\int_U e({\bv})\cdot e({\bv}) d{\bx} dt\right)^{1/2}.
\end{equation}
This norm is induced by the norm
\begin{equation}
\label{norm2}
\parallel {\bv} \parallel_{H_0^1(U)}\equiv \left(\int_U e({\bv})\cdot e({\bv}) d{\bx} \right)^{1/2}.
\end{equation}
By Korn inequality (\cite{OShY}), (\ref{norm2}) is a norm equivalent to the standard one.
Also, ${\mathcal W}_T$ is dense in ${\mathcal V}$. This can be proved in the same way as e.g. 
Thm. 2.1 in \cite{LM}.

A weak solution $\hat{\bm}^\vep$ of (\ref{aux3}), (\ref{aux4}) is an element of ${\mathcal W}_T$ satisfying 
\begin{equation}
\label{aux3}
\int_0^T\int_U \left(\boldsymbol A^\vep e(\hat{\bm}^\vep)-\boldsymbol B^\vep e(\hat{\bm}_t^\vep)\right)\cdot e({\boldsymbol \phi})d{\bx}dt=\langle {\fbf}-{\bg}^\vep, {\boldsymbol \phi}\rangle_{{\mathcal V},{\mathcal V}^\star} 
\end{equation}
for all ${\boldsymbol \phi}\in {\mathcal V}$. Here,
$$
{\bg}^\vep=
-{\rm div}\left( \boldsymbol A^\vep e(\psi{\bm}^\vep_T)-
\boldsymbol B^\vep e(\psi_t{\bm}_T^\vep)\right).
$$

Equation (\ref{aux3})
can be stated as 
\begin{equation}
\label{aux4}
G^\vep\hat{\bm}^\vep
= {\fbf}-{\bg}^\vep,
\end{equation}
with the operator $G^\vep: {\mathcal W}_T \to {\mathcal V}^\star$. We consider it as an unbounded
operator on ${\mathcal V}$ with the domain ${\mathcal W}_T$. The corresponding bilinear form
is defined as
\begin{equation}
\label{bilinear}
\langle G^\vep{\bu}, {\bv}\rangle=
\int_0^T\int_U \left(\boldsymbol A^\vep e({\bu})-\boldsymbol B^\vep e({\bu}_t)\right)\cdot e({\bv})d{\bx}dt
\end{equation}
for each ${\bu}\in {\mathcal W}_T$, ${\bv}\in {\mathcal V}$. 
Finally, we note that the adjoint operator $G^{\vep, \star}$ with the domain 
$$
{\mathcal W}_0=\{ {\bv}\in {\mathcal W}: {\bv}(0)=0\}
$$
is defined by 
\begin{equation}
\label{bilinear}
\langle G^{\vep, \star}{\bu}, {\bv}\rangle=
\int_0^T\int_U \left(\boldsymbol A^\vep e({\bu})+\boldsymbol B^\vep e({\bu}_t)\right)\cdot e({\bv})d{\bx}dt
\end{equation}

\begin{proposition}
\label{gvep-properties}
\noindent
(i) $G^\vep$ is strongly coercive: 
\begin{equation}
\label{co1}
\langle G^\vep{\bu}, {\bu}\rangle\geq \alpha_1\parallel{\bu}\parallel^2_{\mathcal V}
\end{equation}
where $\alpha_1$ is a constant from (\ref{ellipticity}) (and thus independent of $\vep$);\newline
\noindent
(ii) $G^\vep$ has a bounded inverse satisfying
\begin{equation}
\label{inv1}
\parallel (G^\vep)^{-1}{\fbf}\parallel_{{\mathcal V}}\leq \frac{1}{\alpha_1} \parallel {\fbf}\parallel_{{\mathcal V}^\star}
\end{equation}
for each ${\fbf}\in {\mathcal V}^\star$. \newline
\end{proposition}
\noindent
{\it Proof}. If ${\bu}$ is sufficiently smooth then after integrating by parts in (\ref{bilinear})
we would have
\begin{equation*}
\langle G^\vep{\bu}, {\bu}\rangle=
\int_0^T\int_U \boldsymbol A^\vep e({\bu})\cdot e({\bu})d{\bx}dt
+\frac 12\left(\int_U\boldsymbol B^\vep e({\bu})\cdot e({\bu})d{\bx}\right)(0),
\end{equation*}
where we took into account that ${\bu}(T)=0$. However, for an arbitrary ${\bu}\in {\mathcal W}_T$
the second term in the right hand side may not be well
defined. To bypass this difficulty, observe that for almost all $t\in I_T$,  
$$
\frac 12\left(\int_U\boldsymbol B^\vep e({\bu})\cdot e({\bu})d{\bx}\right)(t)
$$
is finite. For such $t$ we have
\begin{eqnarray*}
&&\int_t^T\int_U \left(\boldsymbol A^\vep e({\bu})\cdot e({\bu})
-\boldsymbol B^\vep e({\bu}_t)\cdot e({\bu}\right)d{\bx}dt\\
&=&
\int_t^T\int_U \boldsymbol A^\vep e({\bu})\cdot e({\bu})d{\bx}dt
+\frac 12\left(\int_U\boldsymbol B^\vep e({\bu})\cdot e({\bu})d{\bx}\right)(t)\\
& \geq& 
\int_t^T\int_U \boldsymbol A^\vep e({\bu})\cdot e({\bu})d{\bx}dt
\geq
\alpha_1 \int_t^T\int_U e({\bu})\cdot e({\bu})d{\bx}dt
\end{eqnarray*}
The last inequality follows from (\ref{ellipticity}).
Using absolute continuity in $t$ of the first and last terms in the above inequality,  we can pass to the limit $t\to 0^+$ and obtain 
\begin{eqnarray*}
&\langle G^\vep {\bu}, {\bu}\rangle =& \lim_{t\to0+}
\int_t^T\int_U\left(\boldsymbol A^\vep e({\bu})\cdot e({\bu})
-\boldsymbol B^\vep e({\bu_t})\cdot e({\bu}\right)d{\bx}dt\\
& \hspace*{1.4cm} \geq &
\lim_{t\to0+} \alpha_1 \int_t^T\int_U e({\bu})\cdot e({\bu})d{\bx}dt= 
\alpha_1\parallel {\bu}\parallel_{\mathcal V}^2, 
\end{eqnarray*}
which proves (i).

(ii) follows from (i). This is known (see, e.g. \cite{jko1}, Lemma 1). We only sketch the proof
for completeness. Since
$G^\vep$ is closed, passing to the limit in (\ref{co1}) we obtain that the image of $G^\vep$ is closed
in ${\mathcal V}^\star$. If this image does not contain all ${\mathcal V}^\star$, then,
because of density of ${\mathcal W}_T$ in ${\mathcal V}$, there is
${\bg}\in {\mathcal V}$ such that $\langle G^\vep{\bu}, {\bg}\rangle=0$ for all ${\bu}\in {\mathcal W}_T$. This yields ${\bg}\in {\mathcal W}_0$ (domain of $G^{\vep, \star}$) and
$G^{\vep, \star}{\bg}=0$. Next we observe that $G^{\vep, \star}$ satisfies (\ref{co1}) which yields
${\bg}=0$ and gives a contradiction. Thus $G^\vep$ is onto. The estimate (\ref{inv1}) follows from (\ref{co1}).

\hspace{14cm}$\blacksquare$

{\it Remark}. (ii) implies existence of the pressure $P_3^\vep\in L^2(I_T, L^2(U))$. This follows using standard arguments from \cite{Te} combined with the inclusion
${\boldsymbol A}^\vep e({\bf m}^{pq})+{\boldsymbol B}^\vep e({\bf m}_t^{pq})\in
L^2(I_T, L^2(U))$. Moreover, $P_3^\vep$ is bounded in $L^2(I_T, L^2(U))$ independent of $\vep$. Therefore, extracting a subsequence if necessary, we can assume that $P_3^{\vep}\to \overline{P}_3$ weakly in $L^2(I_T, L^2(U))$.

\begin{definition}
We say that the sequence $G^\vep$ {\it $G$-converges} to an operator $G: D(G)\subset {\mathcal W}_T\to {\mathcal V}^\star$, if for each ${\fbf}\in {\mathcal V}^\star$ the sequence
${\bu}^\vep=\left(G^\vep\right)^{-1}{\fbf}$ converges to some ${\bu}\in D(G)$ weakly in
${\mathcal W}_T$. In this case we define $G{\bu}={\fbf}$.
\end{definition}
\begin{proposition}
\label{g-limit}
The sequence $G^\vep$ contains a $G$-convergent subsequence. The limiting operator $G$ has the following properties:\newline
\noindent
(i) 
\begin{equation}
\label{co2}
\langle G{\bu}, {\bu}\rangle\geq \alpha_1\parallel{\bu}\parallel^2_{\mathcal V}
\end{equation}
for each ${\bu}\in D(G)$; \newline
\noindent
(ii) $D(G)={\mathcal W}_T$.\newline
\end{proposition}
\noindent
{\it Proof}. Let us write $G^\vep=A^\vep-B^\vep\partial_t$ where
$A^\vep, B^\vep: {\mathcal V}\to {\mathcal V}^\star$ are operators induced by the bilinear forms
$$
a^\vep({\bu}, {\bv})\equiv \int_0^T\int_U {\boldsymbol A}^\vep e({\bu})\cdot e({\bv}) d{\bx} dt,
$$
$$
b^\vep({\bu}, {\bv})\equiv \int_0^T\int_U {\boldsymbol B}^\vep e({\bu})\cdot e({\bv}) d{\bx} dt,
$$
respectively.
Ellipticity assumptions (\ref{ellipticity}) imply that $A^\vep, B^\vep$ are coercive and bounded
with coercivity constants and bounds independent of $\vep$.  In particular coercivity
of $B^\vep$ implies that there exists a bounded inverse $\left(B^\vep\right)^{-1}$ defined on ${\mathcal V}^\star$, satisfying
$\parallel \left(B^\vep\right)^{-1}\parallel\leq \frac{1}{\beta_1}$, where $\beta_1$ is the lower bound from (\ref{ellipticity}).  Therefore, if
$(A^\vep-B^\vep\partial_t){\bu}^\vep={\fbf}$ then $\partial_t{\bu}^\vep=\left(B^\vep\right)^{-1}\left(A^\vep{\bu}-{\fbf}\right)$. 
This implies
$$
\parallel \partial_t{\bu}^\vep\parallel_{\mathcal V}\leq \frac{\alpha_2 }{\beta_1} 
\parallel\left( G^\vep\right)^{-1}{\fbf}\parallel_{\mathcal V}+
\parallel{\fbf}\parallel_{{\mathcal V}^\star}\leq
(\frac{\alpha_2}{\alpha_1\beta_1}  +1)\parallel{\fbf}\parallel_{{\mathcal V}^\star}.
$$
Thus, if $\bu^\vep\to \bu$ weakly in ${\mathcal W}_T$ then
\begin{equation}
\label{m-est}
\parallel {\bu}\parallel_{\mathcal W}\leq
C(\alpha_1, \alpha_2, \beta_1)\parallel{\fbf}\parallel_{{\mathcal V}^\star}.
\end{equation}
Since ${\mathcal V}^\star$ is separable, we can use diagonal procedure to find a subsequence,
non relabeled, such that ${\bu}^\vep=\left(G^\vep\right)^{-1}{\fbf}$ converges weakly in
${\mathcal W}_T$ to ${\bu}=G^{-1}{\fbf}$ for all ${\fbf}$ in a dense subset of ${\mathcal V}^\star$. Inequality
(\ref{m-est}) implies that convergence also holds for all ${\fbf}\in {\mathcal V}^\star$, and the operator $G^{-1}$ is bounded. 

Next, consider a sequence ${\bu}^\vep$ such that $G^\vep{\bu}^\vep={\fbf}$. Then, by the preceding, ${\bu}^\vep$ converges weakly to ${\bu}$, and $G{\bu}={\fbf}$ by definition of
the $G$-limit.
Since $\langle G^\vep{\bu}^\vep, {\bu}^\vep\rangle=\langle {\fbf}, {\bu}^\vep\rangle=
\langle G{\bu}, {\bu}^\vep\rangle$,
we can pass to the limit and obtain 
$$
\lim_{\vep\to 0} \langle G^\vep{\bu}^\vep, {\bu}^\vep\rangle=
\langle G{\bu}, {\bu}\rangle.
$$
This, together with lower semicontinuity of the norm with respect to weak convergence, allows passage to the limit in (\ref{co1}) which yields (\ref{co2}).

To prove (ii), observe first that by (\ref{m-est}), $D(G)\subset {\mathcal W}_T$. 
To prove equality, we first show that $G^{-1}$ is injective: $G^{-1}{\fbf}=0$ implies ${\fbf}=0$.
Arguing by contradictions, suppose that there is ${\bg}\in {\mathcal V}^\star, {\bg}\ne 0$ such
that $G^{-1}{\bg}=0$. Consider the sequence ${\bu}^\vep_g=\left(G^\vep\right)^{-1}{\bg}$.
By definition of $G$, ${\bu}^\vep_g$ converges to zero weakly in ${\mathcal W}_T$. 
Then by (\ref{co1}),  
$$
\langle {\bg}, {\bu}_g^\vep\rangle =\langle G^\vep{\bu}_g^\vep, {\bu}_g^\vep\rangle\geq \alpha_1 \parallel {\bu}^\vep_g\parallel^2_{\mathcal V}.
$$
passing to the limit $\vep\to 0$ we obtain that $\parallel {\bu}^\vep_g\parallel_{\mathcal V}\to 0$
and thus ${\bu}^\vep_g$ converges to zero strongly in ${\mathcal V}$. Next, we use ${\bu}^\vep_{g, t}$ as the test function in the weak formulation of $G^\vep{\bu}^\vep_g={\bg}$. Integrating by parts (which can be justified as in the proof of Proposition \ref{gvep-properties}) and using coercivity
of $B^\vep$ we obtain
$$
\left|\langle {\bg}, {\bu}_{g, t}^\vep\rangle\right|=\left|\langle G^\vep{\bu}_g^\vep, {\bu}_{g, t}^\vep\rangle\right|\geq \beta_1 \parallel {\bu}^\vep_{g, t}\parallel^2_{\mathcal V},
$$
and thus ${\bu}^\vep_{g, t}$ converges to zero strongly in ${\mathcal V}$. Next,
using uniform boundedness of $A^\vep, B^\vep$, we write
$$
\parallel {\bg}\parallel_{{\mathcal V}^\star}=
\parallel A^\vep{\bu}^\vep_g-B^\vep {\bu}^\vep_{g, t}\parallel_{{\mathcal V}^\star}
\leq
\alpha_2 \parallel {\bu}^\vep_g \parallel_{\mathcal V}+
\beta_2 \parallel {\bu}^\vep_{g, t} \parallel_{\mathcal V}
$$
where $\alpha_2, \beta_2$ are constants from (\ref{ellipticity}). Passing to the limit
in the above we deduce ${\bg}=0$, which contradicts assumption ${\bg}\ne 0$.
Thus $G^{-1}$ is injective.

Next we show that $D(G)$ (equivalently, the range of $G^{-1}$) is dense in ${\mathcal V}$. If this were false, there would be
${\bh}\ne 0, {\bh} \in {\mathcal V}^\star$ such that
$\langle {\bh}, G^{-1}{\fbf}\rangle=0$ for all ${\fbf}\in {\mathcal V}^\star$. Let
${\bu}_h=G^{-1}{\bh}$. Choosing
${\fbf}={\bh}$, we obtain using (\ref{co2}):
$$
0=\langle {\bh}, G^{-1}{\bh}\rangle=\langle G{\bu}_h, {\bu}_h\rangle\geq
\alpha_1\parallel {\bu}_h\parallel_{\mathcal V}.
$$
Therefore, ${\bu}_h=0$. Then ${\bh}=0$ by injectivity of $G^{-1}$. This contradicts
the assumption ${\bh}\ne 0$.

Finally, observe that
the norm of ${\mathcal W}_T$ induces a scalar product
$$
({\bu}, {\bv})_{{\mathcal W}_T}=
\int_0^T\int_U\left( e({\bu}_t)\cdot e({\bv}_t)+e({\bu})\cdot e({\bv})\right) d{\bx} dt.
$$
In search of a contradiction, suppose that $D(G)$ is a proper subset of ${\mathcal W}_T$. Then
there is $\overline{\bu}\ne 0, \overline{\bu}\in {\mathcal W}_T$ such that $(\overline{\bu}, G^{-1}{\fbf})_{{\mathcal W}_T}=0$ for all ${\fbf}\in {\mathcal V}^\star$.  Expression
$(\overline{\bu}, {\bv})_{{\mathcal W}_T}$ defines a bounded linear functional $l_{\overline{\bu}}({\bv})$
on $D(G)$ which by Hahn-Banach theorem can be extended to a bounded linear functional
$L_{\overline{\bu}}$ on ${\mathcal V}$ and this extension has same norm as $l_{\overline{\bu}}({\bv})$. Therefore, $(\overline{\bu}, G^{-1}{\fbf})_{{\mathcal W}_T}=0$ implies
$\langle L_{\overline{\bu}}, G^{-1}{\fbf}\rangle=0$ for all ${\fbf}\in {\mathcal V}^\star$.
Density of the range of $G^{-1}$ implies $L_{\overline{\bu}}=0$. But then
$(\overline{\bu}, \overline{\bu})_{{\mathcal W}_T}=0$ which contradicts the assumption
$\overline{\bu}\ne 0$. Thus (ii) is proved.

\hspace{14cm}$\blacksquare$

\begin{proposition}
\label{rhs-prop}
There exists $\fbf\in {\mathcal V}^\star$ such that the sequence of solutions ${\bm}^\vep$ of (\ref{aux1}), with
this choice of the right hand side, contains a subsequence (not relabeled) satisfying\newline
\noindent
(i) ${\bm}^\vep\to 0$ weakly in ${\mathcal W}_T$, \newline
\noindent
(ii) ${\bm}^\vep\to 0$ strongly in $L^2(I_T, L^2(U))$.
\end{proposition}

\noindent
{\it Proof}. By Proposition \ref{ab-rhs-choice}, proved below in Sect. \ref{sect:a-prob2}, we can assume that
 ${\bf m}^\vep_T$ converges to zero weakly in $H_0^1(U)$. 
Then $\psi{\bf m}^\vep_T\to 0$ 
and
$\partial_t\psi{\bf m}^\vep_T\to 0$ weakly in ${\mathcal V}$.

Since
${\bm}^\vep=\hat{\bm}^\vep+\psi{\bf m}^\vep_T$, (i) will be proved if we show that
there is choice of $\fbf$ such that $\hat{\bm}^\vep\to 0$ weakly in
${\mathcal W}_T$. 

Consider (\ref{aux4}). In view of (\ref{ellipticity}) and uniform bounds on ${\bm}^\vep_T$, the sequence
${\bg}^\vep$ is bounded in ${\mathcal V}^\star$. Therefore, the sequence
$\left(G^\vep\right)^{-1}\bg^\vep$ is bounded in ${\mathcal W}_T$, and we can extract a 
subsequence that converges weakly to some $\bq\in {\mathcal W}_T$. By Proposition
\ref{g-limit}, (ii), $\bq\in D(G)$. Therefore, we can choose 
\begin{equation}
\label{choice}
\fbf=G\bq.
\end{equation}
Then
$$
\hat{\bm}^\vep=\left(G^\vep\right)^{-1}(\fbf-\bg^\vep)=
\left(G^\vep\right)^{-1}G\bq-\left(G^\vep\right)^{-1}\bg^\vep.
$$
By Proposition \ref{g-limit}, 
$\left(G^\vep\right)^{-1}G\bq \to \bq$ weakly in ${\mathcal W}_T$, up to extraction of a subsequence.
Hence, $\hat{\bm}^\vep\to 0$. Thus (i) is proved.

To prove (ii), observe first that $\bm^\vep_T$ converges to zero strongly in $L^2(U)$ and thus
$\psi{\bf m}^\vep_T\to 0$ and
$\partial_t\psi{\bf m}^\vep_T\to 0$ strongly in $L^2(I_T, L^2(U))$. Therefore, to prove strong convergence
of $\bm^\vep$ it is enough to prove strong convergence of $\hat\bm^\vep$. 

Next, note that
$\hat\bm^\vep_t$ is bounded in ${\mathcal V}$ independent of $\vep$. Now strong convergence
of $\hat\bm^\vep$ is deduced from (i) and J. L. Lions' compactness theorem (see e.g. \cite{Te}, thm. 2.1, ch.III).

\hspace{14cm}$\blacksquare$

\subsection{Auxiliary problems for $\bn^{pq, \vep}, \bm_T^{pq, \vep}$}
\label{sect:a-prob2}
In this section, $p, q$ are once again fixed, so we drop them to simplify notations, and write ${\bn}^\vep$ instead
of ${\bn}^{pq, \vep}$ and so on.

We seek $\bn^\vep$, $\bm^\vep_T$ satisfying, respectively, the auxiliary problems
\begin{equation}
\label{auxn1}
-{\rm div}\left(
\bA^\vep \left({\boldsymbol I}^{pq}+e({\bn^\vep})\right)\right)-\nabla P^\vep_1=\fbf_1,
\end{equation}
\begin{equation}
\label{auxm1}
-{\rm div}\left(
\boldsymbol B^\vep\left( \boldsymbol I^{pq}+e\left({\bn^\vep}\right)\right)-
e\left(\bm^\vep_T\right)
\right)-\nabla P^\vep_2=\fbf_2,
\end{equation}
with suitably chosen $\fbf_1, \fbf_2$. First, we find $\bn^\vep$ from (\ref{auxn1}). Then this 
$\bn^\vep$ should be plugged into (\ref{auxm1}), and then $\bm^\vep_T$ can be found.
The goal, as before, to choose $\fbf_1, \fbf_2$ so that
$\bn^\vep, \bm^\vep_T$ would have subsequences that converge weakly to zero.
Let
$$
V=\{\bv\in H_0^1: {\rm div}\;\bv=0\},
$$
equipped with the norm (\ref{norm2}).
Given $\fbf_1\in V^\star$, $\bn^\vep\in V$ is a weak solution of (\ref{auxn1}) provided
\begin{equation}
\label{auxn2}
\int_U \bA^\vep \left({\boldsymbol I}^{pq}+e({\bn^\vep})\right)\cdot e(\bw)d\bx=
\langle \fbf_1, \bw\rangle_{V, V^\star}
\end{equation}
for all $\bw\in V$. Weak solutions of (\ref{auxm1}) are defined similarly. This identity can be written as 
an operator equation
\begin{equation}
\label{auxn3}
A^\vep \bn^\vep=\fbf_1-\bg^\vep_1,
\end{equation}
where 
$$
\bg^\vep_1=-{\rm div}\left(\bA^\vep\bI^{pq}\right)\in V^\star,
$$
and $A^\vep: V\to V^\star$ is the operator induced by the bilinear form
$$
a^\vep(\bu, \bv)=\int_U \bA^\vep e(\bu)\cdot e(\bv)d\bx.
$$
Similarly, (\ref{auxm1}) can be written as
\begin{equation}
\label{auxm3}
B^\vep \bm^\vep_T=\fbf_2-\bg^\vep_2,
\end{equation}
where the operator $B^\vep$ is induced by 
the form 
$$
b^\vep(\bu, \bv)=\int_U \bB^\vep e(\bu)\cdot e(\bv)d\bx,
$$
and
$$
\bg^\vep_2=-{\rm div}\left( \bB^\vep(\bI^{pq}+e(\bn^\vep))\right).
$$
By (\ref{ellipticity}), operators $A^\vep, B^\vep$ satisfy
\begin{eqnarray}
&&\langle A^\vep\bu, \bu\rangle_{V, V^\star}\geq \alpha_1 \parallel \bu\parallel^2_V\label{co1-1}\\
&&\parallel A^\vep\bu \parallel_{V^\star}\leq \alpha_2 \parallel \bu\parallel_V\label{co1-2}\\
&&\langle B^\vep\bu, \bu\rangle_{V, V^\star}\geq \beta_1 \parallel \bu\parallel^2_V\label{co2-1}\\
&&\parallel B^\vep\bu \parallel_{V^\star}\leq \beta_2 \parallel \bu\parallel_V\label{co2-2}
\end{eqnarray}
Lax-Milgram lemma implies existence of unique solutions of (\ref{auxn3}), (\ref{auxm3}). These solutions satisfy
$$
\bn^\vep=\left(A^\vep\right)^{-1}(\fbf_1-\bg_1^\vep), \;\;\;\;\; \parallel \bn^\vep\parallel_V\leq \frac{1}{\alpha_1}
\parallel \fbf_1-\bg_1^\vep\parallel_{V^\star},
$$
$$
\bm^\vep_T=\left(B^\vep\right)^{-1}(\fbf_2-\bg_2^\vep), \;\;\;\;\; \parallel \bm^\vep_T\parallel_V\leq \frac{1}{\beta_1}
\parallel \fbf_2-\bg_2^\vep\parallel_{V^\star}.
$$
{\it Remark}.  Existence of the pressures $P_1^{\vep} \in L^2(U)$, $P_2^{\vep} \in L^2(U)$ follows using standard arguments from \cite{Te}.  Moreover, these pressures are bounded in $L^2(U)$ independent of $\vep$. Therefore, extracting a subsequence if necessary, we can assume that $P_j^{\vep}\to \overline{P}_j, j=1,2$ weakly in $L^2(U)$.

\begin{definition}
The sequence of operators $A^\vep:\; V\to V^\star$ G-converges to an operator $A$ if $\left(A^\vep\right)^{-1}{\fbf}$
converges to some $\bu\in V$ weakly in $V$, for each $\fbf\in V^\star$. We also define $\fbf=A\bu$.
\end{definition}
\begin{proposition}
\label{ab-gconv}
The sequences $A^\vep, B^\vep$ contain $G$-convergent subsequences.
\end{proposition} 

\noindent
{\it Proof}. This is known \cite{MT}, thm. 2.
\begin{proposition}
\label{ab-rhs-choice}
There exists $\fbf_1\in V^\star$ (respectively $\fbf_2$) such that, up to extraction of a subsequence, $\bn^\vep$
(respectively $\bm^\vep_T$) converge 
to zero weakly in $V$. 
\end{proposition}

\noindent
{\it Proof.} Consider (\ref{auxn3}). Since $\bg_2^\vep$ is bounded in $V^\star$, and $\left(A^\vep\right)^{-1}$ is bounded independent of $\vep$, the sequence $\left(A^\vep\right)^{-1}\bg_1^\vep$ is bounded in $V$. Thus we can
extract a subsequence that converges weakly in $V$ to some $\overline{\bu}_1\in V$. Choose
\begin{equation}
\label{ch1}
\fbf_1=A\overline{\bu}_1.
\end{equation}
Then 
$$
\bn^\vep=\left(A^\vep\right)^{-1}A\overline{\bu}_1-\left(A^\vep\right)^{-1}\bg_1.
$$
By definition of $A$, the first term in the right hand side converges to $\overline{\bu}_1$ weakly in $V$, and so does the second.
Hence $\bn^\vep\to 0$ weakly in $V$. For (\ref{auxm3}) the procedure is the same. Up to extraction of a subsequence,
$\left(B^\vep\right)^{-1}\bg_2^\vep\to \overline{\bu}_2$ weakly in $V$, and we choose
\begin{equation}
\label{ch2}
\fbf_2=B\overline{\bu}_2.
\end{equation}

\hspace{14cm}$\blacksquare$

\section{Inertial terms in the momentum balance equation}
\label{sect:in-terms}
In the remainder of the paper, we assume that the oscillating
test functions $\bw^\vep$ are defined as follows.
\begin{definition}
\label{def:otf}
Let $\bw\in C^\infty(I_T\times U), {\rm div}\bw=0$ be arbitrary, 
and define $\bw^\vep$ by
(\ref{appr-otf}). In (\ref{appr-otf}), choose $\bn^{pq, \vep}, \bm^{pq, \vep}_T$ that solve, respectively
(\ref{auxn1}), (\ref{auxm1}) with the right hand sides chosen according to
(\ref{ch1}), (\ref{ch2}). Also, let $\bm^{pq, \vep}$ satisfy (\ref{aux1}) with
the right hand side chosen according to (\ref{choice}).
\end{definition}
\begin{proposition}
\label{otf-convergence}
The sequence $\bw^\vep$ defined as above satisfies
\begin{equation}
\label{13}
{\bw}^\vep\to {\bw}~~{\rm in}~L^2(I_T, L^2(U)),~~~~
{\bw}^\vep_t\to {\bw}_t~~{\rm in}~L^2(I_T, L^2(U)).
\end{equation}
Also, $\bw^\vep\in L^\infty(I_T, H_0^1(U))$, and
\begin{equation}
\label{sup-bound}
\parallel \bw^\vep\parallel_{L^\infty(I_T, H_0^1(U))}\leq T^{1/2}
\parallel \bw^\vep_t\parallel_{L^2(I_T, H_0^1(U))}\leq
CT^{1/2}
\end{equation}
with $C$ independent of $\vep$.
\end{proposition}

\noindent
{\it Proof}.  
First we need a formula for the time derivative of $\bw^\vep$. After taking time derivative of (\ref{appr-otf}), integrating by parts in the time convolution (which involves putting time differentiation on $\bw_t$ instead of $\bm^{pq, \vep}$) and using $\bw(T)=0$ we obtain 
\begin{eqnarray}
\label{otf-td}
{\bw}^\vep_t& = & {\bw}_t+
{\bn}^{pq, \vep}({\bx})e({\bw}_t)_{pq}+\\
&& \int_t^T \bm^{pq, \vep}(t-\tau+T)e({\bw}_{tt})_{pq}(\tau)d\tau+
\nabla \phi_t^\vep.\nonumber
\end{eqnarray}
To prove strong convergence of $\bw^\vep_t$ we need to prove that all terms in the right hand side of (\ref{otf-td}) converge to zero strongly in $L^2(I_T, L^2(U))$.

\noindent
{\it Step 1. Show that 
$$
{\bn}^{pq, \vep}e({\bw})_{pq}+
\int_t^T \bm^{pq, \vep}(t-\tau+T)e({\bw}_t)_{pq}(\tau)d\tau\to 0,
$$
and
$${\bn}^{pq, \vep}({\bx})e({\bw}_t)_{pq}+
\int_t^T \bm^{pq, \vep}(t-\tau+T)e({\bw}_{tt})_{pq}(\tau)d\tau \to 0$$
strongly in $L^2(I_T\times U)$.}

By Propositions \ref{ab-rhs-choice}, $\bn^{pq,\vep}\to 0$ weakly in $H_0^1(U))$ and thus strongly
in $L^2(U)$. By Proposition \ref{rhs-prop},  $\bm^{pq,\vep}\to 0$ strongly in
$L^2(I_T, L^2(U))$, and we conclude.

\noindent
{\it Step 2. Show that $\nabla\phi^\vep\to 0$ and $\partial_t\nabla\phi^\vep\to 0$ strongly in $L^2(I_T\times U)$.}

First we estimate $e(\nabla\phi^\vep)$. 
Note that
$e(\nabla\phi^\vep)_{ij}=\partial_i\partial_j\phi^\vep$ and 
write
$\partial_i\partial_j\phi^\vep =\partial_i\partial_jE\star (\Delta\phi^\vep)+\partial_i\partial_j K\star (\Delta\phi^\vep)$, where $E=c\frac{1}{|x|}$ is a fundamental solution
of the Laplacian, and $K$ is a harmonic function, which depends only on $U$. To estimate $\partial_i\partial_jE\star (\Delta\phi^\vep)$ we use Calder\'on-Zygmund inequality
(see, e.g. \cite{gt}, Theorem 9.9) with $p=2$ and obtain
$$
\parallel \partial_i\partial_jE\star (\Delta\phi^\vep)\parallel_{L^2(I_T\times U)}\leq
\parallel \Delta\phi^\vep\parallel_{L^2(I_T\times U)}
$$
Since $\partial_i\partial_j K$ is a smooth function, there exists a constant $C(U)$ depending only on $U$
such that
$$
\parallel \partial_i\partial_jK\star (\Delta\phi^\vep)\parallel_{L^2(I_T\times U)}\leq
C(U)\parallel \Delta\phi^\vep\parallel_{L^2(I_T\times U)}
$$
Thus
\begin{equation}
\label{phi5}
\parallel \partial_i\partial_j\phi^\vep\parallel_{L^2(I_T\times U)}\leq
(1+C(U))\parallel \Delta\phi^\vep\parallel_{L^2(I_T\times U)}.
\end{equation}

Combining (\ref{phi5}) with (\ref{phi1}) we find
\begin{equation}
\label{phi2}
\int_{I_T}\int_U |e(\nabla\phi^\vep)|^2 d\bx dt\leq
C(\bw, U)\sum_{p, q=1}^3
\left(\parallel \bn^{pq, \vep}\parallel_{L^2(U)}^2+
\parallel \bm^{pq, \vep}\parallel_{L^2(I_T\times U)}^2\right)
\end{equation}
so $e(\nabla \phi^\vep)$ converges to zero strongly in $L^2(I_T\times U)$. This also implies that the components of the Hessian of $\phi^\vep$ converge to zero strongly in $L^2(I_T\times U)$. 

Finally, the standard a priori estimate for the Neumann problem $\Delta\phi^\vep=f^\vep, f\in L^2(U)$, satisfying
$\nabla\phi^\vep\cdot\bnu=0$ on the boundary, yields
\begin{equation}
\label{neu}
\int_U \nabla\phi^\vep\cdot\nabla\phi^\vep d\bx\leq \parallel f^\vep\parallel_{L^2(U)}
\parallel \phi^\vep\parallel_{L^2(U)}\leq C\parallel f^\vep\parallel_{L^2(U)}
\left(\int_U \nabla\phi^\vep\cdot\nabla\phi^\vep d\bx\right)^{1/2}.
\end{equation}
Here, $C$ is the constant in Poincar\'e inequality.
Poincar\'e inequality applies after we impose the condition
$\int_{\partial U} \phi^\vep dS=0$, standard for Neumann problem. Since $f^\vep$, given by the right hand side of (\ref{phi1}), converges to zero strongly in $L^2(I_T, L^2(U))$, (\ref{neu}) implies that $\nabla\phi^\vep$ converges to zero strongly in $L^2(I_T, L^2(U))$.

Differentiating (\ref{phi1}) in $t$ and integrating by parts as in (\ref{otf-td}) we find
\begin{equation}
\label{phi3}
\Delta \phi_t^\vep=-{\bn}^{pq, \vep}\cdot \nabla e({\bw}_t)_{pq}
-\int_t^T \bm^{pq, \vep}(t-\tau+T)\cdot \nabla e({\bw}_{tt})_{pq}(\tau)d\tau.
\end{equation}
Therefore, arguing as above we have
\begin{equation}
\label{phi4}
\int_{I_T}\int_U |e(\nabla\phi_t^\vep)|^2 d\bx dt\leq
C(\bw, U)\sum_{p, q=1}^3
\left(\parallel \bn^{pq, \vep}\parallel_{L^2(U)}^2+\parallel \bm^{pq, \vep}\parallel _{L^2(U)}^2\right),
\end{equation}
which yields $e(\nabla\phi^\vep_t)\to 0$, and then $\nabla\phi^\vep_t\to$ strongly in $L^2(I_T\times U)$.

\noindent
{\it Step 4. Prove (\ref{sup-bound})}.

Since $\bw^\vep(t)=-\int_t^T \bw^\vep_t(\tau) d\tau$, we obtain for almost all $t\in I_T$
$$
\parallel\bw^\vep\parallel_{H_0^1(U)}(t)\leq
\int_0^T \parallel \bw^\vep_t\parallel_{H_0^1(U)}(\tau) d\tau\leq
T^{1/2}\left(\int_0^T \parallel \bw^\vep_t\parallel^2_{H_0^1(U)}(\tau) d\tau\right)^{1/2},
$$ 
and (\ref{sup-bound}) follows.

\hspace{14cm}$\blacksquare$

Next we will need the following lemma (\cite{lions2}, lemma 5.1).
\begin{lemma}
\label{lemma:lions}
Let $g^n, h^n$ converge weakly to $g, h$, respectively in
$L^{p_1}(0,T; L^{p_2}(\Omega)), L^{q_1}(0,T; L^{q_2}(\Omega))$, where
$1\leq p_1, p_2\leq \infty$,
$$
\frac{1}{p_1}+\frac{1}{q_1}=\frac{1}{p_2}+\frac{1}{q_2}=1.
$$
We assume in addition that 
$\partial_t g^n$ is bounded in $L^{1}(0,T; W^{-m,1}(\Omega))$
for some $m\geq 0$ independent of $n$ and
$$
\parallel h^n-h^n(\cdot+\xi, t)\parallel_{L^{q_1}(0,T; L^{q_2}(\Omega))}\to 0
$$
as $|\xi|\to 0$, uniformly in $n$.

Then $g^nh^n$ converges to $gh$ in the sense of distributions on
$\Omega\times (0,T)$.
\end{lemma} 
This lemma can be used to obtain the effective mass conservation equation and to pass to the limit
in the inertial terms in the momentum equation (\ref{1.3}).

\begin{proposition}
\label{first-in}
Let  $\bw^\vep$ be functions from Definition \ref{def:otf}. Then
\begin{equation}
\label{in0}
\lim_{\vep\to 0}\int_{U} \rho^\vep (0, \bx)\bv_0\cdot\bw^\vep~d\vex =
\int_{U} \overline{\rho}(0, \bx)\; \bv_0\cdot \bw~d\vex ,
\end{equation}
\begin{equation}
\label{in11}
\lim_{\vep\to 0}\int_{I_T\times U} \rho^\vep \bv^\vep \cdot \partial_t\bw^\vep~d\vex dt=
\int_{I_T\times U} \overline{\rho}\; \overline{\bv}\cdot \partial_t\bw~d\vex dt,
\end{equation}
and 
\begin{equation}
\label{in2}
\lim_{\vep\to 0}\int_{I_T\times U} \rho^\vep \bv^\vep\otimes \bv^\vep \cdot \nabla \bw^\vep~d\vex dt=
\int_{I_T\times U} \overline{\rho}\;\overline{\bv}\otimes \overline{\bv} \cdot \nabla \bw~d\vex dt.
\end{equation}
\end{proposition}

\noindent
{\it Proof.} By Lemma \ref{lemma:lions}, $\rho^\vep\bv^\vep\to \overline\rho\;\overline{\bv}$
in the sense of distributions, and thus also weakly in $L^2(I_T\times U)$.
By (\ref{13}), $\bw^\vep, \partial_t\bw^\vep$
converge to respectively $\bw, \partial_t\bw$ strongly in  $L^2(I_T\times U)$.
This permits passage to the limit in the products
and yields (\ref{in11}).

Since 
$\bw^\vep-\bw\in C(I_T, L^2(U))$, 
\begin{eqnarray*}
& \parallel\bw^\vep-\bw\parallel_{L^2(U)}(0) &\leq
\int_0^T \parallel \bw^\vep_t-\bw_t\parallel_{L^2(U)}(\tau) d\tau\\
&& \leq
T^{1/2}\left(\int_0^T \parallel \bw^\vep_t-\bw_t\parallel_{L^2(U)}(\tau) d\tau\right)^{1/2}.
\end{eqnarray*}
Noting that $\bw^\vep_t \to \bw_t $  strongly in $L^2(I_T\times U)$, we obtain
that $(\bw^\vep-\bw)(0)\to 0$ strongly in $L^2(U)$. Since $\rho^\vep(0, \bx)$ converges
weakly-$\star$ in $L^\infty(U)$ to $\overline{\rho}(0, \bx)$, strong convergence of
$\bw^\vep(0, \bx)$ permits passage to the limit in the product $\rho^\vep(0, \bx)\bv_0(\bx)\cdot\bw^\vep(0, \bx)$ and yields (\ref{in0}).

Next, fix $j\in \{1,2,3\}$, pick a function $\eta\in C^\infty_0(I_T\times {U})$,
and insert the test function
$(\bw^\vep-\bw)\eta$ into the weak formulation of the mass balance equation.
This yields
\begin{eqnarray}
&
\int_U \rho^\vep(0, \vex)(\bw^\vep-\bw)\eta(0,\vex)d\vex& -
\int_{I_T}\int_U \rho^\vep \partial_t \left((\bw^\vep-\bw)\eta\right) d\vex dt\label{weak-zero}\\
&& -
\int_{I_T}\int_U \rho^\vep {\bv}^\vep\cdot\nabla \left((\bw^\vep-\bw)\eta\right) d\vex dt=0.
\nonumber
\end{eqnarray}
Strong convergence of $(\bw^\vep-\bw)(0)$ to zero implies
\begin{equation}
\label{int-conv1}
\lim_{\vep\to 0}\int_U \rho^\vep(0, \vex)(\bw^\vep-\bw)\eta(0,\vex)d\vex=0.
\end{equation}
Next, note that
\begin{equation}
\label{int-conv2}
\lim_{\vep\to 0}\int_{I_T}\int_U \rho^\vep \partial_t \left((\bw^\vep-\bw)\eta\right) d\vex dt=0
\end{equation}
because $(\bw^\vep_t-\bw_t) \to 0$  strongly in $L^2(I_T\times U)$, and
$\rho^\vep$ is bounded in $L^\infty(I_T\times U)$ independent of $\vep$.
Now from (\ref{weak-zero}), (\ref{int-conv1}) and (\ref{int-conv2}) we deduce
\begin{equation}
\label{int-conv3}
\lim_{\vep\to 0}
\int_{I_T}\int_U \rho^\vep {\bv}^\vep\cdot\nabla \left((\bw^\vep-\bw)\eta\right) d\vex dt=0.
\end{equation}
Since $\rho^\vep\bv^\vep$ is bounded in $L^2(I_T\times U)$, and $\bw^\vep\to \bw$
strongly in $L^2(I_T\times U)$, (\ref{int-conv3}) implies
\begin{equation}
\label{int-conv4}
\lim_{\vep\to 0}
\int_{I_T}\int_U\eta \rho^\vep {\bv}^\vep\cdot\nabla \left(\bw^\vep-\bw\right) d\vex dt=0.
\end{equation}
Since $\eta\in C^\infty_0(I_T\times U)$ is an arbitray test function, $\rho^\vep {\bv}^\vep\cdot\nabla \left(\bw^\vep-\bw\right)\to 0$ in ${\mathcal D}^\prime(I_T\times U)$. 

Next we claim that $\rho^\vep {\bv}^\vep\cdot\nabla \left(\bw^\vep-\bw\right)$ is bounded in $L^2(I_T, L^{5/6}(U))$ independent of $\vep$. This follows from Sobolev imbedding for $\bv^\vep$ and H\"older inequality.
Application of H\"older inequality yields
$$
\left(\int_U \left|\rho^\vep v_k^\vep\partial_k w^\vep_j\right|^sd\bx\right)(t)\leq
\parallel\rho^\vep\parallel^s_{L^\infty(U)}(t)\left(\int_U \left| v_k^\vep\right|^{sq}\right)^{\frac{1}{q}}(t)
\left(\int_U \left| \partial_k w^\vep\right|^{sq^\prime}\right)^{\frac{1}{q^\prime}}(t).
$$
Here $s, q\geq 1$ and $\frac{1}{q}+\frac{1}{q^\prime}=1$.
Hence,
\begin{eqnarray*}
\int_{I_T}\left(\int_U \left|\rho^\vep v_k^\vep\partial_k w^\vep_j\right|^sd\bx\right)^{\frac{2}{s}}dt
& \leq & 
\parallel\rho^\vep\parallel^2_{L^\infty(I_T\times U)}\int_{I_T}\left(\int_U \left| \partial_k w^\vep_j\right|^{sq^\prime}\right)^{\frac{2}{sq^\prime}}(t)
\left(\int_U \left| v_k^\vep\right|^{sq}d\bx\right)^{\frac{2}{sq}}(t)dt
\nonumber\\
& \leq  &
\parallel\rho^\vep\parallel^2_{L^\infty(I_T\times U)}
\parallel \partial_k w^\vep_j\parallel^{2}_{L^\infty(L^{sq^\prime})}
\int_{I_T}\left(\int_U \left| v_k^\vep\right|^{sq}d\bx\right)^{\frac{2}{sq}}dt
\nonumber
\end{eqnarray*}
Therefore
\begin{equation}
\label{lp-lq}
\parallel\rho^\vep v_k^\vep \partial_k w^\vep_j\parallel_{L^2(I_T, L^s(U))}\leq
\parallel\rho^\vep\parallel_{L^\infty(I_T\times U)}
\parallel \partial_k w^\vep_j\parallel_{L^\infty(L^{sq^\prime})}
\parallel v_k^\vep \parallel_{L^2(I_T, L^{sq}(U))}
\end{equation}
We need to choose $s, q$ so that (i) the right hand side of (\ref{lp-lq}) is finite; and
(ii) $\bv^\vep\in L^2(I_T, L^{s^\prime}(U))$, where $\frac{1}{s}+\frac{1}{s^\prime}=1$.
By Sobolev imbedding, $s^\prime\leq 6$, and therefore
\begin{equation}
\label{s1}
s\geq \frac{6}{5}.
\end{equation}
By  (\ref{sup-bound}), $\parallel \partial_k w^\vep_j\parallel^{\frac{2}{sq^\prime}}_{L^\infty(L^{sq^\prime})}$ is finite if
\begin{equation}
\label{s2}
1\leq sq^\prime \leq 2 \Longleftrightarrow 1-\frac{1}{q}\leq s\leq 2-\frac{2}{q},
\end{equation}
($sq^\prime<2$ are allowed because $U$ is bounded).
Also, by Sobolev imbedding 
\begin{equation}
\label{s3}
1\leq sq \leq 6 \Longleftrightarrow \frac{1}{q}\leq s\leq \frac{6}{q}.
\end{equation}
The solution set of inequalities (\ref{s1})--(\ref{s3}) is a non-empty, convex quadrilateral
in the $1/q-s$-plane. For example, we can choose $s=\frac{6}{5}$ and any 
$q$ satisfying $\frac{1}{5}\leq \frac{1}{q}\leq \frac 25$. 
If, for example, $1/q=\frac 25$, then $sq=3$, $sq^\prime=\frac 95$ and (\ref{lp-lq}) becomes
$$
\parallel\rho^\vep \bv^\vep \cdot \nabla \bw^\vep\parallel_{L^2(I_T, L^{\frac 65}(U))}\leq
C \parallel\rho^\vep\parallel_{L^\infty(I_T\times U)}
\parallel \nabla \bw^\vep\parallel_{L^\infty(I_T, L^{\frac 95}(U))}
\parallel \bv^\vep\parallel_{L^2(I_T, L^{3}(U))}.
$$
Since the right hand side of (\ref{lp-lq}) is bounded independent
of $\vep$, the claim is proved. 

Together with (\ref{int-conv4}), this yields $\rho^\vep
\bv^\vep\cdot\nabla(\bw^\vep-\bw)\to 0$ weakly in $L^2(I_T, L^{6/5}(U))$. Therefore,
the weak limit of $\rho^\vep
\bv^\vep\cdot\nabla \bw^\vep$ is the same as the weak limit of $\rho^\vep
\bv^\vep\cdot\nabla \bw$ in $L^2(I_T, L^{6/5}(U))$. By Lemma \ref{lemma:lions},
$\rho^\vep
\bv^\vep\cdot\nabla \bw \to \overline{\rho}\overline{\bv}\cdot\nabla \bw$. Thus
\begin{equation}
\label{w-conv1}
\rho^\vep
\bv^\vep\cdot\nabla \bw^\vep\to \overline{\rho}\;\overline{\bv}\cdot\nabla \bw
\end{equation}
weakly in $L^2(I_T, L^{6/5}(U))$. The bound on $\partial_t[
\rho^\vep
\bv^\vep\cdot\nabla\bw^\vep]$ in a negative Sobolev space follow from the corresponding bounds
on $\rho^\vep\bv^\vep$ and the fact that $\bw^\vep_t$ is bounded in $L^2(I_T, H_0^1(U))$.

Now application of Lemma \ref{lemma:lions} with $g^\vep=\rho^\vep\bv^\vep\cdot\nabla\bw^\vep$, $h^\vep=\bv^\vep$ yields
\begin{equation}
\label{conv-w}
\lim_{\vep\to 0}\int_{I_T\times U} \rho^\vep \bv^\vep\otimes \bv^\vep \cdot \nabla \bw^\vep\eta~d\vex dt=
\int_{I_T\times U} \overline{\rho}\;\overline{\bv}\otimes \overline{\bv} \cdot \nabla \bw\eta~d\vex dt.
\end{equation}
For each $\eta\in C_0^\infty(I_T\times U)$. From Sobolev imbedding and bounds on $\nabla\bw^\vep$, in the same way
as (\ref{lp-lq}) was analyzed, choosing $1/q=2/5, s=6/5$, we obtain
\begin{equation}
\label{lp-lq1}
\parallel\rho^\vep \bv^\vep \otimes \bv^\vep\cdot \nabla \bw^\vep\parallel_{L^2(I_T, L^{\frac 65}(U))}\leq
C \parallel\rho^\vep\parallel_{L^\infty(I_T\times U)}
\parallel \nabla \bw^\vep\parallel_{L^\infty(I_T, L^2(U))}
\parallel \bv^\vep\otimes \bv^\vep \parallel_{L^2(I_T, L^{3}(U))}
\end{equation}
with $C$ independent of $\vep$. Note that the right hand side is bounded independent of $\vep$.
Therefore, $\rho^\vep \bv^\vep\otimes \bv^\vep \cdot \nabla \bw^\vep$ converges to
$\overline{\rho}\;\overline{\bv}\otimes \overline{\bv} \cdot \nabla \bw$ weakly in $L^2(I_T, L^{\frac 65}(U))$. 
This implies
convergence of integrals of $\rho^\vep \bv^\vep\otimes \bv^\vep \cdot \nabla \bw^\vep$ over subsets
of $I_T\times U$, and in particular (\ref{in2}).

\hspace{14cm}$\blacksquare$
\section{Effective deviatoric stress. Proof of the main theorem}
\label{sect:deviatoric}
\begin{theorem}
\label{thm:effective}
There exist a subsequence, not relabeled, and effective material tensors $\overline{\bA}\in L^2(U), \overline{\bB}\in L^2(U)$ and $\overline{\bC}\in L^2(I_T\times U)$
such that for each $\bw\in C_0^\infty(I_T\times U)$ with ${\rm div}\;\bw=0, \bw(T, \bx)=0$,
\begin{eqnarray*}
&& \int_{I_T}\int_U \overline{\bT}\cdot e(\bw)d\bx dt =
\lim_{\vep\to 0}\int_{I_T}\int_U {\bT}^\vep\cdot e(\bw^\vep)d\bx dt \\
&& =
\int_{I_T}\int_U\left(\overline{\bA}e(\overline{\bu})+
\overline{\bB}e(\overline{\bv})+\int_0^t \overline{\bC}(t-\tau, \cdot)e(\overline{\bv})(\tau, \cdot)\right)\cdot e(\bw)d\bx d\tau
\end{eqnarray*}
as $\vep\to 0$ along this subsequence.
\end{theorem}

\noindent
{\it Proof}. The theorem follows from several propositions.
First, we prove that convergence of inertial terms implies compensated compactness of stress.
\begin{proposition}
\label{prop-comp}
Let $\bw^\vep$ be test functions from Definition \ref{def:otf}. Then
$$
\lim_{\vep\to 0}\int_{I_T}\int_U {\bT}^\vep\cdot e(\bw^\vep)d\bx dt=
\int_{I_T}\int_U \overline{\bT}\cdot e(\bw)d\bx dt.
$$
\end{proposition}

\noindent
{\it Proof of the proposition}. First, use $\bw$ as a test function in (\ref{1.3}) and pass to the limit
$\vep\to 0$. Repeated application of Lemma \ref{lemma:lions} in the inertial terms yields
\begin{eqnarray}
&& -\int_{U} \overline{\rho}(0, \bx)\; \bv_0\cdot \bw~d\vex-
\int_{I_T\times U} \overline{\rho}\; \overline{\bv}\cdot \partial_t\bw~d\vex dt\label{bar1}\\
&&
-\int_{I_T\times U} \overline{\rho}\;
\overline{\bv}\otimes \overline{\bv} \cdot \nabla \bw~d\vex dt
+\int_{I_T\times U} \overline{\bT}\cdot e(\bw) d\bx dt=0\nonumber
\end{eqnarray}
Then insert $\bw^\vep$ into (\ref{1.3}) and pass to the limit $\vep\to 0$. By Proposition \ref{first-in}, the integrals corresponding to the inertial terms will converge to the corresponding integrals of the limiting functions $\bar\rho, \bar{\bv}$. This yields
\begin{eqnarray}
&& -\int_{U} \overline{\rho}(0, \bx)\; \bv_0\cdot \bw~d\vex-
\int_{I_T\times U} \overline{\rho}\; \overline{\bv}\cdot \partial_t\bw~d\vex dt\label{bar2}\\
&&
-\int_{I_T\times U} \overline{\rho}\;
\overline{\bv}\otimes \overline{\bv} \cdot \nabla \bw~d\vex dt
+\lim_{\vep\to 0}\int_{I_T\times U} {\bT}^\vep\cdot e(\bw^\vep) d\bx dt=0\nonumber
\end{eqnarray}

Comparison of (\ref{bar1}) and (\ref{bar2}) finishes the proof.

\hspace{14cm}$\blacksquare$

Next, 
using symmetry of $\bA^\vep, \bB^\vep$ we have
$$
\int_{I_T\times U} {\bT}^\vep\cdot e(\bw^\vep) d\bx dt=
\int_{I_T\times U}e({\bu}^\vep)\cdot\left[\bA^\vep e(\bw^\vep)-\bB^\vep e(\bw_t^\vep)\right]d\bx dt.
$$
Differentiation in (\ref{appr-otf}) yields
\begin{eqnarray}
\boldsymbol A^\vep e\left({\bw}^\vep\right)-\boldsymbol B^\vep e\left({\bw}_t^\vep\right) & =&{\mathcal F}_1^{pq, \vep}e({\bw})_{pq}+{\mathcal F}_2^{pq, \vep}e({\bw}_t)_{pq}+
\int_t^T{\mathcal F}_3^{pq, \vep}(t-\tau+T) e({\bw}_t)_{pq}(\tau)d\tau\label{total-flux}
\\
&&+\boldsymbol A^\vep({\bg}^\vep_1+{\bg}^\vep_2)-\boldsymbol B^\vep(\partial_t {\bg}^\vep_1+\partial_t {\bg}^\vep_2)\nonumber\\
&& +\bA^\vep e(\nabla \phi^\vep)-\bB^\vep e(\nabla \phi_t^\vep)
,
\nonumber 
\end{eqnarray}
where
$$
{\mathcal F}_1^{pq, \vep}({\bx})=\bA^\vep\left({\boldsymbol I}^{pq}+e\left({\bn}^{pq, \vep}\right)\right),~~~~~
{\boldsymbol I}^{pq}=\frac{1}{2} \left({\bf e}_p\otimes {\bf e}_q+{\bf e}_q\otimes {\bf e}_p\right),
$$
$$
{\mathcal F}_2^{pq, \vep}({\bx})=\boldsymbol B^\vep\left( \boldsymbol I^{pq}+e\left({\bn}^{pq, \vep}\right)-
e\left(\bm_T^{pq, \vep}\right)\right),
$$
$$
{\mathcal F}_3^{pq, \vep}(t, {\bx})=
\boldsymbol A^\vep e\left(\bm^{pq, \vep}\right)-
\boldsymbol B^\vep e\left(\bm^{pq, \vep}_t\right),
$$
\begin{equation}
\label{g1}
{\bg}^\vep_1=\frac 12\left(\bn^{pq,. \vep}\otimes \nabla e({\bw})_{pq}+\left(\bn^{pq, \vep}\otimes \nabla e({\bw})_{pq}\right)^\trp\right),
\end{equation}
\begin{equation}
\label{g2}
{\bg}^\vep_2=\frac 12 \int_t^T \left(\bm^{pq, \vep}(t-\tau+T)\otimes \nabla e({\bw})_{pq}(\tau)+\left(\bm^{pq, \vep}(t-\tau+T)\otimes \nabla e({\bw})_{pq}(\tau)\right)^\trp\right)d\tau.
\end{equation}
Next we show that the only terms in (\ref{total-flux}) that contribute to the effective stress are the terms containing
${\mathcal F}^{pq, \vep}_j$.
\begin{proposition}
\label{prop-pert}
Let $\bw^\vep$ be as is Definition \ref{def:otf}. Then
\begin{eqnarray}
\label{pert1}
&& \lim_{\vep\to 0}
\int_{I_T}\int_U
e(\bu^\vep)\cdot
\left[
\boldsymbol A^\vep({\bg}^\vep_1+{\bg}^\vep_2)-\boldsymbol B^\vep(\partial_t {\bg}^\vep_1+\partial_t {\bg}^\vep_2)\right]d\bx dt=0, \\
&&
\lim_{\vep\to 0}\int_{I_T}\int_U
e(\bu^\vep)\cdot
\left[\bA^\vep e(\nabla \phi^\vep)-\bB^\vep e(\nabla \phi_t^\vep)\right]
d\bx dt=0.\label{pert2}
\end{eqnarray}
\end{proposition}

\noindent
{\it Proof of the proposition}. Since $e(\bu^\vep)$ converges to $e(\overline{\bu})$ weakly in $L^2(I_T\times U)$, it is enough
to show that all terms in brackets in (\ref{pert1}), (\ref{pert2}) converge to zero strongly in $L^2(I_T\times U)$. 

\noindent
{\it Step 1. Prove (\ref{pert1}).}

By
Proposition \ref{ab-rhs-choice}, $\bn^{pq, \vep}$ converges to zero strongly in $L^2(I_T\times U)$
Therefore, $\bg^\vep_1, \partial_t\bg^\vep_1$ in (\ref{g1}) converge to zero strongly in $L^2(I_T\times U)$.
By Proposition
\ref{rhs-prop}, $\bm^{pq, \vep}$ converge to zero strongly in $L^2(I_T\times U)$. Hence, $\bg^\vep_2$ in (\ref{g2}) converges to zero strongly in $L^2(I_T\times U)$. 
Next, differentiate $\bg^\vep_2$ in $t$ and integrate by parts in the time convolution exactly as in the proof of Proposition
\ref{otf-convergence}. Then
\begin{eqnarray*}
&&\partial_t\bg^\vep_2=
\frac 12 \int_t^T \left(\bm^{pq, \vep}(t-
\tau+T)\otimes \nabla e({\bw}_{tt})_{pq}(\tau)+\left(\bm^{pq, \vep}(t-\tau+T)\otimes \nabla e({\bw}_{tt})_{pq}(\tau)\right)^\trp\right)d\tau.
\end{eqnarray*}
which converges to zero strongly in $L^2(I_T\times U)$. Next, since $\bA^\vep, \bB^\vep$ are bounded pointwise
independent of $\vep$, we deduce that
$$
\boldsymbol A^\vep({\bg}^\vep_1+{\bg}^\vep_2)-\boldsymbol B^\vep(\partial_t {\bg}^\vep_1+\partial_t {\bg}^\vep_2)
$$
converges to zero strongly in $L^2(I_T\times U)$. 

\noindent
{\it Step 2. Prove (\ref{pert2}).}

From (\ref{phi2}) we have $e(\nabla \phi^\vep)\to 0$, and by (\ref{phi4}), $e(\nabla \phi^\vep_t)\to 0$ strongly
$L^2(I_T\times U)$.
Hence,
$\bA^\vep e(\nabla \phi^\vep)$, $\bB^\vep e(\nabla \phi^\vep_t)$ also converge to zero strongly in $L^2(I_T\times U)$.

\hspace{14cm}$\blacksquare$

By Proposition \ref{prop-pert}, 
$$
\int_{I_T} \int_U \overline{\bT}\cdot e(\bw)d\bx dt=
\lim_{\vep\to 0}\int_{I_T} \int_U \bT^\vep\cdot e(\bw^\vep)
d\bx dt=\lim_{\vep\to 0} I(\bu^\vep, \bw^\vep),
$$
where
\begin{equation}
\label{i1}
I(\bu^\vep, \bw^\vep)=
\int_{I_T} \int_U e(\bu^\vep)\cdot
\left( {\mathcal F}^{pq, \vep}_1 e(\bw)_{pq}+ {\mathcal F}^{pq, \vep}_2 e(\bw_t)_{pq}+
\int_t^T{\mathcal F}_3^{pq, \vep}(t-\tau+T) e({\bw}_t)_{pq}(\tau)d\tau\right)  d\bx dt
\end{equation}
\begin{proposition}
\label{prop:I}
There exist the effective tensors
$\overline{\bA}\in L^2(U), \overline{\bB}\in L^2(U)$ and $\overline{\bC}\in L^2(I_T\times U)$
such that, up to extraction of a subsequence,
\begin{eqnarray}
\label{I-lim}
&\lim_{\vep\to 0}I(\bu^\vep, \bw^\vep)&=
-\langle \overline{\bu}, {\rm div}\left(\overline{\mathcal F}^{pq}_1 e(\bw)_{pq}\right)
+{\rm div}\left(\overline{\mathcal F}^{pq}_2 e(\bw_t)_{pq}\right)\rangle\\
&&
-\langle \bu^\vep, {\rm div}\left(\int_t^T\overline{\mathcal F}_3^{pq}(t-\tau+T) e({\bw}_t)_{pq}(\tau)d\tau
d\bx dt\right)\rangle\nonumber\\
&&
=
\int_{I_T}\int _U
\left(\overline{\bA}e(\overline{\bu})+
\overline{\bB}e(\overline{\bu}_t)+
\int_0^t \overline{\bC}(t-\tau, \cdot)e(\overline{\bu}_t)(\tau, \cdot)d\tau\right)\cdot
e(\bw)d\bx dt,\nonumber
\end{eqnarray}
\end{proposition}

\noindent
{\it Proof of the proposition}.
First, we note that ${\mathcal F}_j^{pq, \vep}, P^{pq, \vep}_j, j=1, 2$ are bounded
in $L^2(U)$ independent of $\vep$. Therefore, extracting weakly convergent subsequences, not relabeled,
we have that
${\mathcal F}_j^{pq, \vep}\rightharpoonup \overline{\mathcal F}_j^{pq}$, 
$P^{pq, \vep}_j \rightharpoonup \overline{P}^{pq}$, and passing to the limit in 
(\ref{auxn1}) and (\ref{auxm1}) we obtain
$$
{\rm div}{\mathcal F}_j^{pq, \vep} -\nabla P^{pq, \vep}_j=
\fbf_j={\rm div}\overline{\mathcal F}_j^{pq} -\nabla \overline{P}^{pq}_j.
$$
and thus
\begin{equation}
\label{flux-lim1}
{\rm div}{\mathcal F}_j^{pq, \vep}=\nabla \left(P^{pq, \vep}_j -\overline{P}^{pq}_j\right)
+{\rm div}\overline{\mathcal F}_j^{pq}, ~~~j=1,2.
\end{equation}
Similarly, ${\mathcal F}^{pq, \vep}_3, P^{pq, \vep}$ are bounded in 
$L^2(I_T\times U)$ independent of $\vep$. Therefore, extracting weakly convergent subsequences as before,
we obtain
\begin{equation}
\label{flux-lim2}
{\rm div}{\mathcal F}_3^{pq, \vep} =\nabla \left(P^{pq, \vep}_3 -\overline{P}^{pq}_3\right)
+{\rm div}\overline{\mathcal F}_3^{pq}.
\end{equation}
From (\ref{flux-lim1}), (\ref{flux-lim2}) we deduce
\begin{eqnarray}
\label{flux-lim3}
&{\rm div}\left({\mathcal F}_1^{pq, \vep}e(\bw)_{pq}\right) =&
{\rm div}\left(\overline{\mathcal F}_1^{pq}e(\bw)_{pq}\right)+ \nabla \left[(P^{pq, \vep}_1 -\overline{P}^{pq}_1)e(\bw)_{pq}\right]\\
&&
-\left(P^{pq, \vep}_1 -\overline{P}^{pq}_1\right)\nabla e(\bw)_{pq}
+\left({\mathcal F}_1^{pq, \vep}-\overline{\mathcal F}_1^{pq}\right)\cdot\nabla e(\bw)_{pq}.\nonumber
\end{eqnarray}
and
\begin{eqnarray}
\label{flux-lim4}
&{\rm div}\left({\mathcal F}_j^{pq, \vep}e(\bw_t)_{pq}\right) =&
{\rm div}\left(\overline{\mathcal F}_j^{pq}e(\bw_t)_{pq}\right)+ 
\nabla \left[(P^{pq, \vep}_j -\overline{P}^{pq}_j)e(\bw_t)_{pq}\right]\\
&&
-\left(P^{pq, \vep}_j -\overline{P}^{pq}_j\right)\nabla e(\bw_t)_{pq}
+\left({\mathcal F}_j^{pq, \vep}-\overline{\mathcal F}_j^{pq}\right)\cdot\nabla e(\bw_t)_{pq},
~~~j=2, 3.\nonumber
\end{eqnarray}
Next, since ${\rm div}\;\bu^\vep=0$,
\begin{eqnarray}
\label{pressure}
&&\langle \bu^\vep, \nabla \left[(P^{pq, \vep}_1 -\overline{P}^{pq}_1)e(\bw)_{pq}\right]\rangle=0,\\
&& \langle \bu^\vep, \nabla \left[(P^{pq, \vep}_2 -\overline{P}^{pq}_2)e(\bw_t)_{pq}\right]\rangle=0,
\nonumber\\
&& \langle \bu^\vep, \nabla \left[\int_t^T (P^{pq, \vep}_3 -\overline{P}^{pq}_3)(t-\tau)
e(\bw_t)_{pq}(\tau)d\tau\right]\rangle=0.\nonumber
\end{eqnarray}

Integrating by parts and using (\ref{flux-lim3}), (\ref{flux-lim4}) and (\ref{pressure}) we have
\begin{eqnarray*}
&I(\bu^\vep, \bw^\vep)
= &-\langle \bu^\vep, {\rm div}\left(\overline{\mathcal F}^{pq}_1 e(\bw)_{pq}\right)
+{\rm div}\left(\overline{\mathcal F}^{pq}_2 e(\bw_t)_{pq}\right)\rangle\\
&&
-\langle \bu^\vep, {\rm div}\left(\int_t^T\overline{\mathcal F}_3^{pq}(t-\tau+T) e({\bw}_t)_{pq}(\tau)d\tau
d\bx dt\right)\rangle
+{\mathcal R}^\vep
\end{eqnarray*}
where
\begin{eqnarray*}
&{\mathcal R}^\vep=&
- \langle\bu^\vep,  ({\mathcal F}^{pq, \vep}_1-\overline{\mathcal F}^{pq}_1)\cdot\nabla e(\bw)_{pq}\rangle
-\langle ({\mathcal F}^{pq, \vep}_2-\overline{\mathcal F}^{pq}_2)\cdot\nabla e(\bw_t)_{pq}\rangle\\
&&
-\langle \bu^\vep, \int_t^T\left({\mathcal F}_3^{pq, \vep}-\overline{\mathcal F}_3^{pq}\right) (t-\tau+T) \cdot \nabla e({\bw}_t)_{pq}(\tau)d\tau
d\bx dt\rangle\\
&&
-\langle\bu^\vep,  (P^{pq, \vep}_1-\overline{P}^{pq}_1)\cdot\nabla e(\bw)_{pq}\rangle
-\langle\bu^\vep,  (P^{pq, \vep}_2-\overline{P}^{pq}_2)\cdot\nabla e(\bw_t)_{pq}\rangle\\
&&
-\langle \bu^\vep, \int_t^T\left(P_3^{pq, \vep}-\overline{P}_3^{pq}\right) (t-\tau) \cdot \nabla e({\bw}_t)_{pq}(\tau)d\tau
d\bx dt\rangle
\end{eqnarray*}
since $\bu^\vep \in L^2(I_T, H_0^1(U))$, we can apply
Lemma \ref{lemma:lions} which yields $\lim_{\vep\to 0}{\mathcal R}^\vep=0$. Then we have
(\ref{I-lim}),
where the components of the effective tensors $\overline{\bA}(\bx), \overline{\bB}(\bx), \overline{\bC}(t, \bx)$ are defined by
\begin{equation}
\label{effective1}
\overline{A}_{pqij}=\overline{\mathcal{F}}^{pq}_{1, ij},
~~~~~~
\overline{B}_{pqij}=\overline{\mathcal{F}}^{pq}_{2, ij},
~~~~~~
\overline{C}_{pqij}=\overline{\mathcal{F}}^{pq}_{3, ij}.
\end{equation}

\hspace{14cm}$\blacksquare$

This completes the proof of the theorem \ref{thm:effective}.

\noindent
{\it Proof of the main theorem}.

To obtain (\ref{main2}), we pass to the limit in (\ref{1.2}) using Lemma \ref{lemma:lions}. Next, insert
$\bw^\vep=N^\vep\bw$ into (\ref{1.3}). The limit of the inertial terms is given in Proposition \ref{first-in}, and the limit
of the term containing $\bT^\vep\cdot e(\bw^\vep)$ is provided by Theorem \ref{thm:effective}. Together, these results yield
(\ref{main3}). The divergence-free constraint (\ref{main1}) is obtained by straightforward passing to the limit
in ${\rm div}\;\bv^\vep=0$.

\hspace{14cm}$\blacksquare$
\section{Fluid-structure interaction}
\label{sect:fsi}
Compared to the previous sections, the main difference now is lack of ellipticity in $\bA^\vep$. In this section we assume
$\bA^\vep=\bA_1 \theta_0^\vep$. This means that phase one is a Kelvin-Voight viscoelastic material, and phase two is a Newtonian fluid. To deal with degeneration of $\bA^\vep$
we modify (\ref{appr-otf}) as follows.
\begin{eqnarray}
\label{appr-otf-f}
N^\vep\bw\equiv {\bw}^\vep& = & {\bw}+\int_t^T \bn^{pq, \vep}(t-\tau+T)e({\bw})_{pq}(\tau)d\tau+\\
&&
\int_t^T \bm^{pq, \vep}(t-\tau+T)e({\bw}_t)_{pq}(\tau)d\tau+ \nabla \phi^\vep.\nonumber
\end{eqnarray}
Here, $\bn^{pq, \vep}, \bm^{pq, \vep}, \phi^\vep$ are as in (\ref{appr-otf}), (\ref{phi1}), respectively. We note, 
however, that the auxiliary problems for $\bn^{pq, \vep}, \bm^{pq, \vep}$ will be different.
Differentiating (\ref{appr-otf-f}) we obtain
\begin{eqnarray}
\boldsymbol A^\vep e\left({\bw}^\vep\right)-\boldsymbol B^\vep e\left({\bw}_t^\vep\right) & =&{\mathcal F}_1^{pq, \vep}e({\bw})_{pq}+
\left({\mathcal F}_2^{pq, \vep}\right) ({\bw}_t)_{pq}+\label{total-flux1}
\\
&&\int_t^T{\mathcal F}_3^{pq, \vep}(t-\tau+T) e({\bw})_{pq}(\tau)d\tau+
\int_t^T{\mathcal F}_4^{pq, \vep}(t-\tau+T) e({\bw}_t)_{pq}(\tau)d\tau
\nonumber\\
&&+\boldsymbol A^\vep({\bg}^\vep_1+{\bg}^\vep_2)-\boldsymbol B^\vep(\partial_t {\bg}^\vep_1+\partial_t {\bg}^\vep_2)\nonumber\\
&& +\bA^\vep e(\nabla \phi^\vep)-\bB^\vep e(\nabla \phi_t^\vep),
\nonumber 
\end{eqnarray}
where
$$
{\mathcal F}_1^{pq, \vep}=\bA^\vep{\boldsymbol I}^{pq}-\bB^\vep e\left({\bn}^{pq, \vep}_T\right),~~~~~
{\boldsymbol I}^{pq}=\frac{1}{2} \left({\bf e}_p\otimes {\bf e}_q+{\bf e}_q\otimes {\bf e}_p\right),
$$
$$
{\mathcal F}_2^{pq, \vep}=-\bB^\vep\left( \boldsymbol I^{pq}+e\left({\bm}^{pq, \vep}_T\right)
\right),
$$

$$
{\mathcal F}_3^{pq, \vep}({\bx})=\boldsymbol A^\vep e\left(\bn^{pq, \vep}\right)-
\boldsymbol B^\vep e\left(\bn^{pq, \vep}_t\right),
$$

$$
{\mathcal F}_4^{pq, \vep}(t, {\bx})=
\boldsymbol A^\vep e\left(\bm^{pq, \vep}\right)-
\boldsymbol B^\vep e\left(\bm^{pq, \vep}_t\right),
$$
\begin{equation}
\label{g1}
{\bg}^\vep_1=\frac 12 \int_t^T \left(\bn^{pq, \vep}(t-\tau+T)\otimes \nabla e({\bw})_{pq}(\tau)+\left(\bn^{pq, \vep}(t-\tau+T)\otimes \nabla e({\bw})_{pq}(\tau)\right)^\trp\right)d\tau.
\end{equation}
\begin{equation}
\label{g2}
{\bg}^\vep_2=\frac 12 \int_t^T \left(\bm^{pq, \vep}(t-\tau+T)\otimes \nabla e({\bw}_t)_{pq}(\tau)+\left(\bm^{pq, \vep}(t-\tau+T)\otimes \nabla e({\bw}_t)_{pq}(\tau)\right)^\trp\right)d\tau.
\end{equation}
In the statements of the following auxiliary problems we drop $p, q$ to simplify notations.
\begin{itemize}
\item{\it First auxiliary problem}. 
Find $\bn^\vep_T \in H_0^1(U)$ satisfying
\begin{equation}
\label{f-aux1}
{\rm div}\left( \bA^\vep \bI^{pq}-\bB^\vep e\left({\bn}^{\vep}_T\right) \right)-\nabla P^{\vep}_1=\fbf_1,\;\;\;\;\;\;{\rm div}\;\bn^\vep_T=0,
\end{equation}
\item{\it Second auxiliary problem}. 
Find $\bm^\vep_T \in H_0^1(U)$ satisfying
\begin{equation}
\label{f-aux2}
-{\rm div}\left( \bB^\vep (\bI^{pq}+e\left({\bm}^{\vep}_T\right))\right)-\nabla P^{\vep}_2=\fbf_2,\;\;\;\;\;\;{\rm div}\;\bm^\vep_T=0,
\end{equation}
\item{\it Third auxiliary problem}. 
Find $\bn^\vep\in {\mathcal W}$ satisfying
\begin{equation}
\label{f-aux3}
-{\rm div}\left( (\bA^\vep-\bB^\vep\partial_t) e\left({\bn}^{\vep}\right)\right)-\nabla P^{\vep}_3=\fbf_3,\;\;\;\;\;\;{\rm div}\;\bn^\vep=0,\;\;\;\;\;\;\bn^\vep(T)=\bn^\vep_T.
\end{equation}
\item{\it Fourth auxiliary problem}. Find $\bm^\vep\in {\mathcal W}$ satisfying
\begin{equation}
\label{f-aux4}
-{\rm div}\left( (\bA^\vep-\bB^\vep\partial_t) e\left({\bm}^{\vep}\right)\right)-\nabla P^{\vep}_4=\fbf_4,\;\;\;\;\;\;{\rm div}\;\bm^\vep=0,\;\;\;\;\;\;\bm^\vep(T)=\bm^\vep_T.
\end{equation}
\end{itemize}

Since $\bB^\vep$ is still elliptic, the problems (\ref{f-aux1}), (\ref{f-aux2}) can be analyzed exactly as problems in
Section \ref{sect:a-prob2}. All the results in that section apply without change. The problems (\ref{f-aux3}), (\ref{f-aux4}) 
were dealt with in Section \ref{sect:a-prob1}. The most important condition is still ellipticity of $\bB^\vep$, and the only change that is needed is in the proof of (i) in Proposition \ref{gvep-properties} where we used
ellipticity of $\bA^\vep$. Now, to prove (i) note that
\begin{equation}
\label{i1}
\int_t^T\int_U\bA^\vep e(\bu)\cdot e(\bu)d \bx dt+\frac 12 \int_U\bB^\vep e(\bu)\cdot e(\bu)d \bx(t)=\int_t^T\langle \fbf, \bu\rangle_{H^{-1}(U), H_0^1(U)}(\tau)d\tau
\end{equation}
holds for almost all $t\in I_T$. Estimating the right hand side of (\ref{i1}) we have
\begin{eqnarray}
\label{i2}
&&\int_t^T\langle \fbf, \bu\rangle_{H^{-1}(U), H_0^1(U)}(\tau)d\tau\leq
\int_t^T\parallel\fbf\parallel_{H^{-1}(U)}(\tau)\parallel\bu\parallel_{H_0^1(U)}(\tau)d\tau\\
&& \leq
\left(\int_t^T\parallel\fbf\parallel^2_{H^{-1}(U)}(\tau)d \tau\right)^{1/2}
\left(\int_t^T\parallel\bu\parallel^2_{H_0^1(U)}(\tau)d \tau\right)^{1/2}\nonumber\\
&&
\leq
\parallel\fbf\parallel_{L^2(I_T, H^{-1}(U))}\parallel\bu\parallel_{L^2(I_T, H_0^1(U))}\nonumber
\end{eqnarray}
Since the first term in the left hand side of (\ref{i1}) is non-negative, 
from (\ref{i1}), (\ref{i2}) we deduce using ellipticity of $\bB^\vep$
\begin{equation*}
\frac 12 \beta_1\sup_{t\in I_T}\int_U e(\bu)\cdot e(\bu)d \bx(t)\leq
\parallel\fbf\parallel_{L^2(I_T, H^{-1}(U))}\parallel\bu\parallel_{L^2(I_T, H_0^1(U))}.
\end{equation*}
This yields (i) (with a different constant $\frac 12 \beta_1 T^{-1}$) after observing that
$$
\parallel \bu\parallel^2_{L^2(I_T, H_0^1(U))}\equiv 
\int_0^T\int_U e(\bu)\cdot e(\bu)d \bx\leq T \sup_{t\in I_T}\int_U e(\bu)\cdot e(\bu)d \bx(t).
$$

All the arguments in Sections \ref{sect:in-terms}, \ref{sect:deviatoric} apply with minor changes due to the presence of
four fluxes ${\mathcal F}_j$, $j=1,\ldots, 4$, instead of three.  The result can be summarized as a theorem.

\begin{theorem}
\label{fsi-main}
In the case if fluid-structure interaction,
the limits
$\overline{\rho}, \overline{\bv}$, $\overline{\bu}$ satisfy 
\begin{equation}
\label{main1}
{\rm div}\;\overline{\bv}=0,
\end{equation}
and the integral identities
\begin{equation}
\label{main2}
\int_U \overline{\rho}_0\phi(0,\vex)d\vex-\int_{I_T}\int_U \overline{\rho}\partial_t \phi d\vex dt-
\int_{I_T}\int_U \overline{\rho}\;\overline{\bv}\cdot\nabla \phi d\vex dt=0,
\end{equation}
\begin{eqnarray}
&&-\int_U \overline{\rho}_0{\bv}_0 \cdot {\boldsymbol \psi} d\vex-\int_{I_T}\int_U \overline{\rho}\;
\overline{\bv}\cdot\partial_t {\boldsymbol
\psi} d\vex dt
-\int_{I_T}\int_U \overline{\rho}\;\overline{\bv}\otimes \overline{\bv}\cdot \nabla{\boldsymbol \psi}d\vex dt
\label{main3}
\\
&&+\int_{I_T}\int_U \overline{\bT}\cdot e({\boldsymbol \psi})d\vex dt
=0.\nonumber
\end{eqnarray}
for all smooth test functions 
$\phi$, ${\boldsymbol \psi}$, such that ${\rm div}\;\bpsi=0$, and $\phi, \bpsi$ are equal to zero
on $\partial U$ and vanish for $t\geq T$. 

Moreover, there exist the effective tensors
$\overline{\bA}\in L^2(U), \overline{\bB}\in L^2(U)$ and $\overline{\bC}\in L^2(I_T\times U), \overline{\bD}\in L^2(I_T\times U)$ such that the effective deviatoric stress $\overline{\bT}$ satisfies
 
\begin{equation}
\label{main4}
\overline{\bT}=
\overline{\bA}e(\overline{\bu})+
\overline{\bB}e(\overline{\bv})+
\int_0^t \overline{\bC}(t-\tau)e(\overline{\bu})(\tau)d\tau+
\int_0^t \overline{\bD}(t-\tau)e(\overline{\bv})(\tau)d\tau
\end{equation}
\end{theorem}

\noindent
{\it Remark}. The effective tensors are obtained as weak $L^2(I_T\times U)$ limits of the four fluxes:
\begin{equation}
\label{effective1}
\overline{A}_{pqij}=\overline{\mathcal{F}}^{pq}_{1, ij},
~~~~~~
\overline{B}_{pqij}=\overline{\mathcal{F}}^{pq}_{2, ij},
~~~~~~
\overline{C}_{pqij}=\overline{\mathcal{F}}^{pq}_{3, ij},
~~~~~~
\overline{D}_{pqij}=\overline{\mathcal{F}}^{pq}_{3, ij},
\end{equation}

\section{Acknowledgments}
Work of Alexander Panchenko was supported in part by DOE grant DE-FG02-05ER25709 and by NSF grant
DISE-0438765.

\begin{appendix}
\section{Moving interface in the inertial terms and frozen interface
in the constitutive equations}
\label{ap-model}
In this section we present formal calculations leading to the weak formulation
of the momentum balance equation.

\noindent
{\bf 1. Inertial terms.} 
\begin{eqnarray}
&& \int_{I_T}\int_{V^\vep}\partial_t (\rho^\vep{\bv}^\vep)\cdot \bpsi d{\bx}dt=
\int_{I_T}\int_{U}\theta^\vep \partial_t(\rho^\vep{\bv}^\vep)\cdot \bpsi d{\bx}dt=\label{a1}\\
&&-\int_{I_T}\int_{U}\theta^\vep (\rho^\vep{\bv}^\vep)\cdot\partial_t\bpsi d{\bx}dt-
\int_{I_T}\int_{U}(\rho^\vep{\bv}^\vep\cdot\bpsi)\partial_t\theta^\vep  d{\bx}dt-\nonumber\\
&& \int_{U}\rho^1\theta_0{\bv}_0\cdot\bpsi(0, {\bx}) d{\bx}dt.\nonumber
\end{eqnarray}
\begin{eqnarray}
&& \int_{I_T}\int_{V^\vep}\partial_j(\rho^\vep v^\vep_i v^\vep_j)\psi_i d{\bx}dt=
\int_{I_T}\int_{U}\theta^\vep\partial_j(\rho^\vep v^\vep_i v^\vep_j)\psi_i d{\bx}dt=\label{a2}\\
&&-\int_{I_T}\int_{U}\theta^\vep \rho^\vep v^\vep_i v^\vep_j
\partial_j \psi_i d{\bx}dt
-\int_{I_T}\int_{U}(\rho^\vep {\bv}^\vep\cdot\bpsi) ({\bv}^\vep\cdot
\nabla \theta^\vep) d{\bx}dt.\nonumber
\end{eqnarray}
Combining (\ref{a1}) and (\ref{a2}) we obtain
\begin{eqnarray}
&& \int_{I_T}\int_{V^\vep}\left[\partial_t (\rho^\vep{\bv}^\vep)
+{\rm div}(\rho^\vep{\bv}\otimes{\bv})\right]\cdot \bpsi d{\bx}dt=
\label{a3}\\
&&-\int_{U}\rho^1\theta_0{\bv}_0\cdot\bpsi(0, {\bx}) d{\bx}dt
-\int_{I_T}\int_{U}\theta^\vep (\rho^\vep{\bv}^\vep)\cdot\partial_t
\bpsi d{\bx}dt
-\int_{I_T}\int_{U}\theta^\vep (\rho^\vep{\bv}^\vep)\cdot
\partial_t\bpsi d{\bx}dt\nonumber\\
&&
\int_{I_T}\int_{U}(\rho^\vep{\bv}^\vep\cdot\bpsi)(\partial_t\theta^\vep 
+{\bv}^\vep\cdot\nabla \theta^\vep) d{\bx}dt.
\nonumber
\end{eqnarray}
When $\theta^\vep$ satisfies the interface evolution equation, the last term
in the right hand side is zero. If the interface were frozen, then this term
would be present in the weak formulation of the momentum equation. Since
$\theta_0$ is piecewise constant, the weak formulation would contain a non-physical term supported on the interface. To avoid such non-physical terms, one
needs to use $\theta^\vep$ in the inertial terms of the momentum equation.

\noindent
{\bf 2. Constitutive equation}.
Moving interface assumption combined with the Hook's law would
lead to a non-physical dissipation of the elastic energy.
Indeed, let the elastic part of the stress be written as
\begin{equation}
\label{ac-const0}
\theta_0^\vep(\bx){\boldsymbol A}^1e({\bu}^\vep)+
(1-\theta_0^\vep(\bx)){\boldsymbol A}^2e({\bu}^\vep),
\end{equation}
The important
condition here is that $\theta_0$ is independent of $t$. The stiffness tensors ${\boldsymbol A}^1$,
${\boldsymbol A}^2$ of the phases are supposed to be constant.
Formally multiplying (\ref{ac-const0}) by ${\bv}^\vep$ 
and integrating by parts we obtain
\begin{eqnarray}
\int_0^T\int_U (\theta_0 {\boldsymbol A}^1+(1-\theta_0){\boldsymbol A}^2) 
e({\bu}^\vep)\cdot e(\partial_t {\bu}^\vep)d{\bx}dt\label{f-bal}\\
=
\frac 12\int_0^T\int_U \partial_t\left[(\theta_0 {\boldsymbol A}^1+(1-\theta_0){\boldsymbol A}^2) 
e({\bu}^\vep)\cdot e({\bu}^\vep)\right]d{\bx}dt\nonumber\\
=
\frac 12\int_U\left.(\theta_0 {\boldsymbol A}^1+(1-\theta_0){\boldsymbol A}^2) 
e({\bu}^\vep)\cdot e({\bu}^\vep)d{\bx}\right\arrowvert_0^T.\label{cons}
\end{eqnarray} 
This expresses the fact that the elastic energy changes by the amount of work
done by elastic forces, with no dissipation. If one were to use $\theta^\vep(t, {\bx})$ in (\ref{f-bal}), differentiation
in time would not commute with multiplication by $\theta^\vep {\boldsymbol A}^1+(1-\theta^\vep){\boldsymbol A}^2$, and  (\ref{cons}) could not be obtained. 
\section{Existence of weak solutions. Outline of the proof.}
\label{ap-exist}
In this Section we outline the proof of existence of global weak solutions
for the system (\ref{1.2}), (\ref{1.3}) for each fixed $\vep>0$. Since $\vep$ is fixed, we
drop superscript to simplify the notation. We follow closely 
\cite{lions2}, Sect. 2.3, 2.4. 

1. The initial conditions.  The initial conditions satisfy (\ref{in-den}), (\ref{in-vel}).

2. Formal a priori estimates (\ref{ap-den})--(\ref{ap-flux1}) are obtained
as explained in Sect. \ref{weak-sect}. In particular, renormalization as in \cite{lions2}, 
Sect 2.3 is used to get 
$|\{{\bx}\in U: \alpha\leq \rho (t, {\bx})\leq \beta\}|$ for each $0\leq \alpha\leq \beta<\infty$, where $|\cdot|$ denotes Lebesgue measure. This implies
(\ref{ap-den}).

3. Compactness results.
Since we wish to approximate exact solutions, a compactness result is needed. Suppose
that we have two sequences $\rho^n, {\bv}^n$ satisfying the conditions $0\leq \rho^n\leq C$\;\;${\rm div}~{\bv}^n=0$, a. e. on $I_T\times U$, $\parallel {\bv}^n\parallel_{L^2(I_T, H^1_0(U))}\leq C$, ${\bv}^n \rightharpoonup {\bv}$ weakly in $L^2(I_T, H^1_0(U))$. Moreover, assume that 
$$
\partial_t \rho^n +{\rm div}(\rho^n{\bv}^n)=0
$$
in ${\mathcal D}^\prime (I_T\times U)$, $\rho^n|{\bv}^n|^2$ is bounded in $L^\infty(I_T, L^1(U)$, and we have
$$
\left|\langle \partial_t(\rho^n {\bv}^n), \phi\rangle\right|\leq C \parallel \phi\parallel_{L^q(I_T, W^{m,q}(U))}
$$
for all $\phi\in L^q(I_T, W^{m,q}(U))$ such that ${\rm div}~\phi=0$.

Finally suppose that the initial conditions for the density $\rho^n_0$ satisfy
$$
\rho^n_0 \to \rho_0
$$
in $L^1(U)$.

Then, by the compactness Theorem 2.4 in \cite{lions2} (it applies without any change),
we have 
\begin{eqnarray}
&& \rho^n \to \rho ~~~{\rm in}~C([0,T], L^p(U)),~{\rm for}~{\rm all}~1\leq p <\infty,
\nonumber\\
&& \sqrt{\rho^n}v^n_i \to \sqrt{\rho}v_i
~~~{\rm in}~L^p(I_T, L^r(U)),~{\rm for}~2<p <\infty, 1\leq r\leq \frac{6p}{3p-4},
\nonumber\\
&& v^n_i \to v_i
~~~{\rm in}~L^\theta(I_T, L^{3\theta}(U)),~{\rm for}~1\leq \theta < 2
~~~~{\rm on}~{\rm the}~{\rm set}~\{\rho >0\},
\nonumber
\end{eqnarray}

4. Construction of smooth approximate solutions. As in \cite{lions2}, Sect. 2.4, 
construct solutions $(\rho, {\bv})$ of the approximate system
\begin{eqnarray}
&& \partial_t \rho +{\rm div}(\rho{\bv}_\delta)=0\label{sm-density}\\
&&
\partial_t (\rho {\bv})+{\rm div}(\rho {\bv}_\delta \otimes {\bv})-
{\rm div}\left({\boldsymbol A}_\delta e({\bu})+{\boldsymbol B}_\delta
e({\bv})\right)+\nabla P_\delta=0~~{\rm in}~{\mathcal D}^\prime\label{sm-moment}\\
&& {\rm div}~{\bv}=0~~{\rm in}~{\mathcal D}^\prime, \nonumber
\end{eqnarray}
where ${\bv}_\delta$, ${\boldsymbol A}_\delta$, ${\boldsymbol B}_\delta$, are smooth
regularizations of the respective quantities. The initial conditions are regularizations of the original ones. Then, using a fixed point argument as in
Theorem 2.6, we can prove existence of smooth solutions to (\ref{sm-density}), 
(\ref{sm-moment}). Existence of a fixed point follows from the a priori estimates.
The only issue that needs to be explained here is bootstrap regularity of the
constructed solutions. The following proposition replaces Proposition 2.1 in \cite{lions2}
\begin{proposition}
\label{reg-prop}
Consider
the system
\begin{equation}
\label{reg1}
c\partial_t v_i+b \cdot\nabla v_i-a\Delta v_i-m\Delta u_i+\partial_i P=0,
\end{equation}
${\rm div}~{\bv}=0$ in $I_T\times U$, $i=1,2,3$, with the initial conditions
\begin{equation}
\label{reg2}
{\bv}(0,\cdot)={\bv}_0,~~~~~{\bu}(0, \cdot)=0.
\end{equation}
Suppose that $c\in L^\infty(I_T\times U)$, $a, m\in L^\infty(I_T, W^{1,\infty}(U))$,
$b\in L^2(I_T, L^\infty(U))$, $c\geq k$, $a\geq k$, $m\geq k$ a.e. on $I_T\times U$
for some $k>0$; ${\bv}_0\in H_0^1(U)$. Also, assume that $a, m$ are independent of $t$.

Then the system (\ref{reg1}), (\ref{reg2}) has a unique solution $({\bv}, P)$
such that ${\bv}\in L^2(I_T, H^2(U))\cap C([0,T], H_0^1(U)$, 
$\partial_t {\bv}\in L^2(I_T\times U)$, $\nabla P\in L^2(I_T\times U)$.
\end{proposition}

\noindent
{\it Outline of the proof}. Compared to the proof of Proposition 2.1 in \cite{lions2}, we have one new term $m\Delta u_i$. Multiplying this term by $\partial_t v_i$ and integrating by parts we have
$$
-\int_U m\Delta u_i \partial_t v_i=d_t \int_U \nabla u_i\cdot m\nabla v_i-
\int_U m|\nabla v_i|^2+\int_U \nabla u_i\cdot \nabla m \partial_t v_i 
$$
Multiplying (\ref{reg1}) by $\partial_t v_i$ using the above identity, integrating by parts and summing over $i$ we find
$$
\int_U c|\partial_t{\bv}|^2+\frac 12d_t\int_U a|\nabla{\bv}|^2+
d_t\int_U\nabla {\bu}\cdot m\nabla {\bv}-\int_U m|\nabla {\bv}|^2=
\int_U {\fbf},
$$
where ${\fbf}=b\cdot\nabla{\bv}\partial_t {\bv}-\nabla {\bu}\cdot (\nabla m\otimes \partial_t {\bv})$.
Integrating this identity with respect to $t$ we obtain
\begin{eqnarray}
&& k \int_0^t\int_U |\partial_t {\bv}|^2 +
k \frac 12\int _U |\nabla{\bv}|^2(t)\nonumber \\
&& \leq
\int_U m|\nabla {\bu}||\nabla {\bv}|(t)+\int_0^t\int _U m|\nabla {\bv}|^2+
\int_0^t\int_U {\fbf}+\frac 12\int_U a|\nabla {\bv}|(0).\nonumber
\end{eqnarray}
Next we write  
$$
\int_U m|\nabla {\bu}||\nabla {\bv}|(t)\leq \parallel m\parallel_{L^\infty(U)}
\left( \frac 12 \nu \int_U |\nabla {\bv}|^2(t)+\frac{1}{2\nu} \int_U |\nabla {\bu}|^2(t) \right),
$$
where we choose $\nu=k/2$. The term containing ${\fbf}$ is handled similarly, putting 
$\nu$ in front of $\int_0^t\int_U |\partial_t {\bv}|^2$. Combining the previous two inequalities with
the standard a priori bounds on $\nabla {\bu}, \nabla {\bv}$ we have
\begin{eqnarray}
&& \frac{k}{4} \int_0^t\int_U |\partial_t {\bv}|^2 +
\frac{k}{4}\int _U |\nabla{\bv}|^2(t)
\leq
C_{k},\nonumber
\end{eqnarray}
where  $C_{k}$ depends only on $k$ and the data. This yields a priori estimates on
${\bv}$ in $L^\infty(I_T, H_0^1(U))$ and on $\partial_t {\bv}$ in
$L^2(I_T\times U)$.  

Now we can write (\ref{reg1}) as
\begin{equation}
\label{reg11}
\Delta(a{\bv}+m{\bu})-\nabla P=h,~~~~{\rm div}~{\bv}=0,~~~{\rm in}~U,
\end{equation}
${\bv}\in H_0^1(U), {\bu}\in H_0^1(U)$ for almost all $t\in I_T$. Also,
$h$ is bounded in $L^2(I_t\times U)$ in terms of the data. Next, estimating pressure $P$ exactly as in
\cite{lions2}, Prop. 2.1, we conclude that
$a{\bv}+m{\bu}\in L^2(I_T, H^2(U))$. Since $a, m$ are smooth and positive,
this implies
\begin{equation}
\label{reg3} 
\partial_t {\bu}+\frac{m}{a}{\bu}={\bg},
\end{equation}
where ${\bg}\in L^2(I_T, H^2(U))$. Now, formally,
${\bu}(t,\cdot)=\int_0^t e^{-\frac{m}{a}(t-\tau)}{\bg}(\tau, \cdot)d\tau$, and hence
$$
\left|\partial^2_{x_i x_j}{\bu}\right|^2(t)
\leq \left(\int_0^t \left|\partial^2_{x_i x_j}{\bg}\right|(\tau)d\tau\right)^2
\leq t\int_0^t \left|\partial^2_{x_i x_j}{\bg}\right|^2(\tau)d\tau
\leq T\int_0^t \left|\partial^2_{x_i x_j}{\bg}\right|^2(\tau)d\tau.
$$
Integrating over $U$ we deduce that ${\bu}$ is bounded in $L^\infty(I_T, H^2(U))$
in terms of the data. Then using ${\bv}=-\frac{m}{a}{\bu}+{\bg}$ we obtain
that ${\bv}$ is bounded in $L^2(I_T, H^2(U))$ in terms of the data.
The formal calculations can be easily justified by an approximation argument.
 
\hspace{14cm}$\blacksquare$

5. Passage to the limit. This is done using compactness from step 3 exactly as in
\cite{lions2}.

\end{appendix}


\end{document}